\documentclass[11pt]{article}
\usepackage{amsmath,amsbsy,amsfonts,amssymb}

\usepackage[frenchb]{babel}

\begin{document}

\title{Classification et Changement de base pour les  s\'eries discr\`etes des groupes unitaires 
p-adiques}
\author{C. M{\oe}glin}
\date{}
\maketitle
Il est bien clair pour l'auteur que les m\'ethodes qu'Arthur d\'eveloppent pour les groupes classiques s'appliquent plut\^ot plus simplement pour les groupes unitaires et cet article n'est donc pas original et il utilise les indications donn\'ees par Arthur. Son \'ecriture se justifie pour disposer plus rapidement de cas o\`u l'analyse fine que nous faisons des paquets d'Arthur soit utilisable sans conjectures et pour que les hypoth\`eses de \cite{europe} et \cite{ams} soient d\'emontr\'ees dans le cas des groupes unitaires. D'autre part, les m\'ethodes que nous employons diff\`erent tr\`es certainement par endroits des m\'ethodes d'Arthur;  ce que nous obtenons et qui n'est pas annonc\'e par Arthur (m\^eme si on pourrait le d\'eduire de ce qui est annonc\'e) est le fait que l'on caract\'erise les paquets de Langlands par des propri\'et\'es d'irr\'eductibilit\'e de certaines induites, c'est ce que nous avons appel\'e les blocs de Jordan d'une s\'erie discr\`ete.

Ici on suppose que les lemmes fondamentaux ordinaires ou plus exactement le transfert est valide pour les groupes unitaires et leurs groupes endoscopiques et les groupes lin\'eaires tordus comme expliqu\'e ci-dessous et leurs groupes endoscopiques. Les premi\`eres hypoth\`eses sont disponibles gr\^ace \`a \cite{laumonngo} \cite{waldspurgerfondamental} et \cite{waldspurgertransfert}. Les deuxi\`emes le sont gr\^ace \`a \cite{waldspurgernouveau} et \cite{waldspurgernouveau2}.

Apr\`es ce pr\'eambule, expliquons ce que nous avons en vue  ici: on suppose que $F$ est un corps p-adique et que $E$ est une extension quadratique de $F$ non scind\'ee. On note $\theta$ l'automorphisme ext\'erieur de $GL(n,E)$ qui \`a $g\in GL(n,E)$ associe 
$J ^t\overline{g}^{-1} J^{-1}$, o\`u $J$ est la matrice antidiagonale avec \`a la place $i$ l'\'el\'ement $(-1)^{i+1}$. On dit qu'une repr\'esentation temp\'er\'ee de $GL(n,E)$ est $\theta$-discr\`ete si sa classe de conjugaison est invariante sous l'action de $\theta$ et si elle n'est pas une induite propre \`a partir d'un parabolique $\theta$-stable et d'une repr\'esentation $\theta$-invariante de ce parabolique. On consid\`ere le groupe produit semi-direct de $GL(n,E)$ avec le groupe $\{1,\theta\}$ et on note $\tilde{G}_{n}$ la composante de $\theta$ dans ce groupe. Une th\'eorie de l'endoscopie pour $\tilde{G}_{n}$ a \'et\'e 
d\'evelopp\'ee par Langlands, Kottwitz et Shelstad; on connait les groupes endoscopiques, ou plut\^ot les donn\'ees endoscopiques pour $\tilde{G}_{n}$. Le but de l'article est double: d'abord on montre l'analogue de \cite{arthurnouveau} \sl Local Induction Hypothesis \rm (p. 244) c'est-\`a-dire que pour $\pi$ une repr\'esentation $\theta$-discr\`ete de $GL(n,E)$ prolong\'ee en une repr\'esentation du produit semi-direct, en une repr\'esentation $\tilde{\pi}$,  il existe exactement une 
donn\'ee endoscopique $<H>$ de $\tilde{G}_{n}$ tel que $\tilde{\pi}$ soit le changement de base pour un paquet de repr\'esentations de $H$ suivant cette donn\'ee endoscopique. Cela permet de d\'efinir les paquets stables de s\'eries discr\`etes pour $U(n,E/F)$,  de fa\c{c}on tr\`es concr\`ete, un ensemble de 
s\'eries discr\`etes ${\cal J}$ est un paquet stable, s'il existe une s\'erie $\theta$-discr\`ete $\pi$ de 
$GL(n,E)$ qui soit le changement de base stable d'une combinaison lin\'eaire convenable d'\'el\'ements de ${\cal J}$. On caract\'erise les 
$\pi$ qui peuvent intervenir \`a la Langlands par des homomorphismes de 
$W_{E}\times SL(2,{\mathbb C})$ dans le groupe $GL(n,{\mathbb C})$. Puis
 on montre que l'ensemble des s\'eries discr\`etes de $U(n,E/F)$ est la r\'eunion disjointe des paquets stables et on calcule le nombre d'\'el\'ements dans chaque paquet; cela donne une bijection num\'erique entre les \'el\'ements d'un paquet et les caract\`eres d'un groupe (qui peut s'identifier \`a un centralisateur de $\psi$ dans un groupe convenable). On donne quelques pr\'ecisions  sur cette bijection; en particulier elle v\'erifie les propri\'et\'es de \cite{europe} et \cite{ams}. Au passage on donne une classification des repr\'esentations cuspidales de $U(n,E/F)$ en termes de param\`etres de Langlands. Ceci prouve les hypoth\`eses de base utilis\'ees dans \cite{europe} et \cite{ams}.
 
L'endoscopie tordue a \'et\'e \'etudi\'e par Kottwitz et Shelstad en \cite{ks} et aussi par Labesse  (\cite{labesseJIMJ}. On a repris ici les d\'efinitions de Labesse et en particulier \cite{asian} qui d\'ecrit explicitement les groupes endoscopiques qui interviennent ici. Dans cette r\'ef\'erence, les choix sont 
pr\'ecis\'es. Ces choix nous donne un caract\`ere 
$\omega$ de $E^*$ dont la restriction \`a $F^*$ est le caract\`ere de $F^*$ correspondant \`a l'extension $E$ de $F$. Ainsi le changement de base instable se d\'eduit du changement de base stable par tensorisation par le caract\`ere $\omega$.
Plus g\'en\'eralement, on note $D(n)$ l'ensemble des couples ordonn\'es d'entiers, $n_{1},n_{2}$  dont l'un peut \^etre nul tels que $n=n_{1}+n_{2}$. Pour $(n_{1},n_{2})\in D(n)$, on note $s_{n_{1},n_{2}}$ l'\'el\'ement de $GL(n,{\mathbb C})$ diagonal dont les $n_{1}$ premi\`eres valeurs propres sont $1$ et les $n_{2}$ derni\`eres sont $-1$. Il d\'efinit naturellement un \'el\'ement du groupe dual de $GL(n,E)\rtimes \{1,\theta\}$ et donc une donn\'ee endoscopique (cf. \cite{asian} 1.4); 
le groupe endoscopique correspondant est  $H_{n_{1}, n_{2}}$ le produit des 2 groupes unitaires quasid\'eploy\'es, $U(n_{1},E/F)\times U(n_{2},E/F)$. Cela \'epuise l'ensemble des donn\'ees endoscopiques pour $\tilde{G}_{n}$ (cf. loc.cit.). Si $n_{2}=0$, $H_{n,0}=U(n)$ fournit le changement de base stable tandis que si $n_{1}=0$, $H_{0,n}=U(n)$ fournit le changement de base instable. Le 
r\'esultat de cet article sur le changement de base, est de montrer qu'\`a toute 
repr\'esentation $\theta$-discr\`ete $\pi$ de $GL(n,E)$ correspond exactement une donn\'ee endoscopique $n_{1},n_{2}$ tel qu'il existe $\Pi^{H_{n_{1},n_{2}}}$ un ensemble de repr\'esentations de $H_{n_{1},n_{2}}$ tel que $tr\, \tilde{\pi}$ soit un transfert endoscopique d'une  distribution de la forme 
$\sum_{\tau\in \Pi^{H_{n_{1},n_{2}}}}c_{\tau}tr\, \tau$ o\`u les $c_{\tau}$ sont des nombres complexes.

Si $F$ \'etait le corps des r\'eels, un tel r\'esultat analogue a \'et\'e d\'emontr\'e par Labesse-Harris (\cite{asian}) en s'appuyant sur les r\'esultats de Clozel (\cite{clozelens}); le r\'esultat d\'emontr\'e en loc. cit. se limite au cas o\`u $n_{1}n_{2}=0$. 

Revenons au cas o\`u $F$ est un corps p-adique; on utilise l'espace des int\'egrales orbitales des fonctions cuspidales sur $\tilde{G}_{n}$, not\'e $I_{cusp}(\tilde{G}_{n})$; cette id\'ee revient \`a Arthur \cite{arthurselecta}. On d\'efinit un espace analogue pour les groupes $H_{n_{1},n_{2}}$, not\'e $I_{cusp}(H_{n_{1},n_{2}})$ et un sous-espace de cet espace form\'e par les int\'egrales orbitales stables, $I_{cusp}^{st}(H_{n_{1},n_{2}})$. Il faut recopier \cite{arthurselecta} et \cite{waldspurger1} pour obtenir une 
d\'ecomposition en somme directe:
$$
I_{cusp}(\tilde{G}_{n})=\oplus_{(n_{1},n_{2})\in D(n)} I_{cusp}^{st}(H_{n_{1},n_{2}}) \eqno(1)
$$
l'application de $I_{cusp}(\tilde{G}_{n})$ vers le membre de droite est la somme des transferts (dont on 
a maintenant l'existence);  ici aucun automorphisme ext\'erieur n'appara\^{\i}t dans la situation. En suivant Arthur \cite{arthurselecta} 3.5 on  v\'erifie que cette application est une isom\'etrie pour un produit scalaire convenablement d\'efini en loc. cit. En suivant nos r\'ef\'erences et pr\'ecis\'ement ici \cite{waldspurger1}, on sait associer \`a une 
repr\'esentation temp\'er\'ee $\theta$-discr\`ete de $GL(n,E)$ comme ci-dessus un \'el\'ement de 
$I_{cusp}(\tilde{G}_{n})$; cet \'el\'ement, $f_{\tilde{\pi}}$, est obtenu en consid\'erant la projection sur 
$I_{cusp}$ d'un pseudo coefficient (qui doit \^etre d\'efini) de l'extension $\tilde{\pi}$. De plus 
$I_{cusp}(\tilde{G}_{n})$ est lin\'eairement engendr\'e par ces \'el\'ements $f_{\tilde{\pi}}$. Le r\'esultat que nous avons en vu, s'exprime exactement par le fait que l'image de $f_{\tilde{\pi}}$ dans le membre de droite de (1) est nulle pour tout \'el\'ement de $D(n)$ sauf pr\'ecis\'ement l'un d'entre eux. 

Le r\'esultat qualitatif pr\'ec\'edent se pr\'ecise de la fa\c{c}on suivante permettant de calculer 
l'\'el\'ement $D(n)$ en question: on \'ecrit $\pi$ comme une induite de repr\'esentations de Steinberg 
g\'en\'eralis\'ees:
$$
\pi=\times_{(\rho,a)}St(\rho,a), \eqno(2)
$$
o\`u $(\rho,a)$ parcourt un ensemble de couples o\`u $\rho$ est une repr\'esentation cuspidale 
irr\'eductible  d'un groupe $GL(d_{\rho},E)$, telle que $\rho\simeq \, ^\theta \rho$ o\`u $\theta$ est l'analogue du $\theta$ d\'ej\`a d\'efini pour $n$ remplac\'e par $d_{\rho}$ et o\`u $a$ est un entier. On dit que $(\rho,a)$ est stable si

$a$ est pair et la fonction $L$ d'Asai-Shahidi associ\'ee \`a $\rho$ a un p\^ole en $s=0$

$a$ est impair et la fonction $L$ d'Asai-Shahidi associ\'ee \`a $\rho$ n'a pas de p\^ole en $s=0$.

On entend par fonction $L$ d'Asai-Shahidi la fonction $L$ de \cite{shahidi} qui contr\^ole la
r\'eductibilit\'e de l'induite de $\rho$ vu comme repr\'esentation du parabolique de Levi $GL(n,E)$ du groupe d\'eploy\'e $U(2n,E/F)$. Cette fonction a \'et\'e \'etudi\'ee par Goldberg en \cite{goldberg}  et c'est ''esssentiellement'' la fonction $L$ d'Asai usuelle; le essentiellement cache un point non trivial qui est que la fonction $L$ d'Asai est associ\'ee \`a une repr\'esentation de $W_{F}$ et non $W_{E}$. Il faut donc d\'ej\`a passer de $\rho$ \`a une repr\'esentation de $W_{F}$ et non pas $W_{E}$. Evidemment ici on n'a pas besoin de cette interpr\'etation, une fa\c{c}on \'el\'ementaire de remplacer l'existence d'un p\^ole est de dire que l'induite de $\rho$ \`a $U(2d_{\rho},E/F)$ est 
irr\'eductible (existence d'un p\^ole) ou non irr\'eductible absence de p\^ole. C'est un des r\'esultats de Goldberg qui inclut en particulier que soit l'induite de $\rho$ est irr\'eductible soit c'est l'induite de $\omega\otimes \rho$ qui a cette propri\'et\'e et les 2 options sont exclusives l'une de l'autre (cf. introduction de \cite{goldberg} (1)). Dans le cas o\`u $d_{\rho}$ est impair, la situation est plus simple car contr\^oler par le caract\`ere central de $\rho$ mais  de fa\c{c}on d\'esagr\'eable notre fonction $L$ est en fait une torsion de la fonction $L$ d'Asai usuelle par le caract\`ere quadratique de $W_{F}$ qui correspond \`a l'extension $E$ de $F$. Pour \'eviter les confusions on note $L(\rho,r'_{A},s)$ cette fonction $L$ ind\'ependamment de la parit\'e de $d_{\rho}$.

Avec les notations ci-dessus, on note ${\cal I}_{\pi}^{st}$ l'ensemble des couples $(\rho,a)$ apparaissant dans (2) tels que $(\rho,a)$ soit stable.
Alors $f_{\tilde{\pi}}$ a une projection nulle dans $I_{cusp}^{st}(H_{n_{1},n_{2}})$ sauf exactement si $n_{1}=\sum_{(\rho,a)\in {\cal I}^{st}_{\pi}}ad_{\rho}$.

\

On sait donc maintenant ce qu'est une repr\'esentation elliptique stable de $\tilde{G}_{n}$; elle est associ\'ee \`a un morphisme $\psi$ de $W_{E} \times SL(2,{\mathbb C})$ dans $GL(n,{\mathbb C})$, $\theta$-invariant, sans multiplicit\'e et qui se d\'ecompose en somme de repr\'esentations irr\'eductibles correspondant \`a des s\'eries discr\`etes $St(\rho,a)$ (cf. ci-dessus) telles que $(\rho,a)$ soit stable. On note $\Pi(\psi)$ le paquet de s\'eries discr\`etes de $U(n,E/F)$ associ\'ees \`a $\psi$. On montre que 
$\vert \Pi(\psi)\vert=2^{\ell(\psi)-1}$, o\`u $\ell(\psi)$ est la longueur de $\psi$ vu comme repr\'esentation de $W_{E}\times SL(2,{\mathbb C})$. On montre aussi que $\Pi(\psi)$ ne contient de repr\'esentations cuspidales que si $\psi$ est sans trou; on peut exprimer cette condition de la fa\c{c}on suivante. Soit $a$ un entier, on note $\psi[a]$ la composante isotypique de la repr\'esentation $\psi$ pour la 
repr\'esentation irr\'eductible de $SL(2,{\mathbb C})$ de dimension $a$. Pour tout $a$, $\psi[a]$ est naturellement une repr\'esentation de $W_{E}$. On dit que $\psi$ est sans trou si pour tout $a>2$, $\psi[a]$ est une sous-repr\'esentation de $\psi[a-2]$ en tant que repr\'esentation de $W_{E}$. En termes concrets, cela dit que si $St(\rho,a)$ est une des composantes de $\pi(\psi)$, avec $a>2$, alors $St(\rho,a-2)$ en est aussi une. Pour finir, on associe \`a tout \'el\'ement de $\Pi(\psi)$ un caract\`ere d'un sous-groupe de $(Cent_{GL(n,{\mathbb C})}\psi)^\theta$ (on reprend les id\'ees de \cite{europe}); ce sous-groupe et ce caract\`ere sont canoniquement d\'efinis et le caract\`ere refl\`ete les propri\'et\'es des modules de Jacquet des \'el\'ements de $\Pi(\psi)$. En particulier un tel sous-groupe admet au plus un caract\`ere altern\'e (on ne donne pas la d\'efinition dans cette introduction); si $\psi$ est sans trou, l'existence d'un caract\`ere altern\'e est \'equivalente \`a l'existence de repr\'esentations cuspidales dans $\Pi(\psi)$ puisque les repr\'esentations cuspidales sont alors exactement celles dont le caract\`ere associ\'e est ce caract\`ere altern\'ee. Cela permet de calculer exactement le nombre de repr\'esentations cuspidales dans un paquet $\Pi(\psi)$. Si $\psi$ n'est pas sans trou, cela se g\'en\'eralise en rempla\c{c}ant cuspidale par s\'eries discr\`etes fortement positives (cf. \cite{europe}).

Les conjectures d'Arthur associent \`a tout \'el\'ement de $\Pi(\psi)$ un caract\`ere de tout le groupe; on 
v\'erifie que les conditions mises par Arthur assurent que la restriction de ce caract\`ere au sous-groupe que nous avons d\'efini est notre caract\`ere.

\section{D\'efinitions}
\subsection{Le groupe $GL(n)$ tordu et ses classes de conjugaison stable}
On recopie ici \cite{waldspurger1}.
On reprend les notations $n$, $F$, $E$ de l'introduction ainsi que $\theta$ l'automorphisme ext\'erieur de $GL(n,E)$; on appelle $G^+_{n}$ le produit semi-direct $GL(n,E)\rtimes \{1,\theta\}$ et comme dans l'introduction, on note $\tilde{G}_{n}$ la composante connexe de $\theta$ dans $G^+_{n}$. Soient deux 
\'el\'ements $g,g' \in \tilde{G}_{n}$ semi-simples. On suit la d\'efinition de Labesse \cite{labesseJIMJ} concernant la conjugaison stable et suivant cette d\'efinition $g,g'$ sont stablement conjugu\'es si et seulement si il existe $x\in GL(n, \overline{F})\rtimes \theta$ tel que $g=xg'x^{-1}$ et pour tout \'el\'ement $\sigma\in Gal (\overline{F}/F)$, $\sigma(x^{-1})x$ est dans le sous-groupe du centralisateur de $g$ engendr\'e par la composante connexe de ce stabilisateur et le commutant de $\tilde{G}$ dans 
$GL(n,\overline{F})$. Supposons maintenant que $g$ est fortement r\'egulier, c'est-\`a-dire que son centralisateur est commutatif et la composante neutre de ce centralisateur est un tore. On note 
$\tilde{G}_{n,reg}$ l'ensemble des \'el\'ements fortement r\'eguliers de $\tilde{G}_{n}$.

On dit que $g$ est elliptique si $gZ_{GL(n,E)}(g)^0$ n'est inclus dans aucun sous-groupe de Levi de
$\tilde{G}_{n}$; il faut donc pr\'eciser ce que l'on entend par sous-groupe de Levi de $\tilde{G}_{n}$; ici la d\'efinition est limpide, on consid\`ere les sous-groupes de Levi $\theta$-stables, $M$, de $GL(n,E)$ et on d\'efinit ais\'ement $\tilde{M}$; ces groupes $\tilde{M}$ forment l'ensemble des sous-groupes de Levi de $\tilde{G}_{n}$.  On note $\tilde{G}_{n,ell}$ l'ensemble des \'el\'ements fortement 
r\'eguliers et elliptiques de $\tilde{G}$. Pour toute fonction $f$ localement constante \`a support compact sur $\tilde{G}$, on sait d\'efinir l'int\'egrale orbitale $J^{GL(n,E)}(g,f)$ pour tout $g\in \tilde{G}_{reg}$. On dit que $f$ est une fonction cuspidale si l'int\'egrale orbitale $J^{GL(n,E)}(g,f)=0$ pour tout 
$g\in \tilde{G}_{n,reg}-\tilde{G}_{n,ell}$. On note $C_{cusp}(\tilde{G}_{n})$ l'ensemble des fonctions cuspidales. On note $I_{cusp}(\tilde{G})$ l'ensemble des applications $$I(f): \quad g\in  \tilde{G}_{n,ell} \mapsto J^{GL(n,E)}(g,f)$$pour $f$ parcourant $I_{cusp}(\tilde{G})$. C'est un espace vectoriel.

\subsection{Pseudo-coefficients \label{pseudocoefficients}} 
Pour pouvoir utiliser \cite{waldspurger1} qui \'etend la th\'eorie de Schneider-Stuhler au groupe non connexe, il faut juste changer une d\'efinition, celle de $Z_{2}$ avec la notation de loc. cit. Dans la 
th\'eorie de \cite{ss} il faut que le centre du groupe soit compact et l'id\'ee de \cite{waldspurger1} est donc de ne consid\'erer que des \'el\'ements invariants sous l'action du groupe engendr\'e par 
l'\'el\'ement $z_{2}:=z(\varpi)\theta(z(\varpi))^{-1}$ o\`u ici $\varpi$ est une uniformisante de $E$ et 
$z(\varpi)$ est la matrice diagonale de $GL(n,E)$ de coefficients diagonaux tous \'egaux \`a $\varpi$.  La remarque sous-jacente est que l'on 
s'int\'eresse aux repr\'esentations $\theta$-discr\`ete irr\'eductible de $GL(n,E)$; une telle 
repr\'esentation \`a donc un caract\`ere central invariant par $\theta$ et $z_{2}$ y agit donc trivialement. On note donc ici $Z_{2}$ le sous-groupe de $GL(n,E)$ engendr\'e par $z_{2}$; $Z_{2}$ est invariant sous l'action de $\theta$ et agit par multiplication \`a droite sur $\tilde{G}$. L'autre remarque (\cite{waldspurger1} fin de II.2) est  la multiplication \`a droite par $z_{2}^{-1}$ n'est autre que la conjugaison sous $z(\varpi)$. Comme on ne travaille qu'avec des objets invariants par conjugaison sous $GL(n,E)$, ces objets sont naturellement invariants par multiplication \`a droite par le groupe $Z_{2}$. On peut alors recopier le corollaire \cite{waldspurger1} II.2 en se limitant aux repr\'esentations irr\'eductibles. Dans l'\'enonc\'e ci-dessous $\Delta$ est le d\'eterminant de Weyl usuel.

\

\bf Th\'eor\`eme. \sl Soit $\tilde{\pi}$ une repr\'esentation irr\'eductible de $G^+_{n}$ dont la restriction \`a $GL(n,E)$ reste irr\'eductible. Alors il existe une fonction $f_{\tilde{\pi}}\in C_{cusp}(\tilde{G})$ tel que pour tout \'el\'ement $g\in \tilde{G}_{n,ell}$:
$$
tr\, \tilde{\pi}(g)=\Delta(g)^{-1/2} J^{GL(n,E)}(g,f_{\tilde{\pi}}).
$$\rm

\subsection{Repr\'esentations elliptiques\label{baseelliptique}}
On reprend \cite{waldspurger1} qui a \'et\'e \'ecrit dans un cadre englobant facilement notre cas; c'est ici IV.5 de loc.cite qui nous int\'eresse et ce r\'esultat \'etend des r\'esultats d'Arthur 
\cite{arthuracta}. Pour toute repr\'esentation  $\theta$-discr\`ete $\pi$ de $GL(n,E)$ fixons un prolongement $\tilde{\pi}$ \`a $G^+_{n}$ et une fonction cuspidale $f_{\tilde{\pi}}$ satisfaisant au 
th\'eor\`eme de \ref{pseudocoefficients}. On note $I_{\tilde{\pi}}$ l'image de $f_{\tilde{\pi}}$ dans 
$I_{cusp}(\tilde{G})$.

\

\bf Th\'eor\`eme (Arthur, Waldspurger). \sl L'ensemble des \'el\'ements $I_{\tilde{\pi}}$ quand $\pi$ parcourt l'ensemble des 
repr\'esentations $\theta$-discr\`etes de $GL(n,E)$ forme une base de $I_{cusp}(\tilde{G})$.\rm

\subsection{D\'ecomposition de $I_{cusp}(\tilde{G})$ \label{decomposition}}
En \cite{arthurselecta} 3.5, Arthur a stabilis\'e l'espace $I_{cusp}(G)$ quand $G$ est un groupe r\'eductif connexe d\'eploy\'e et quand on conna\^{\i}t l'existence du transfert pour tous les groupes endoscopiques elliptiques de $G$. Pr\'ecis\'ement, le transfert montre que pour tout $f\in I_{cusp}(G)$ et pour toute 
donn\'ee endoscopique elliptique $<H>$, il existe $f^{<H>}\in I_{cusp}(H)$ tel que $f^{<H>}$ soit un transfert de $f$ et soit stable; on peut en plus imposer \`a $f^{<H>}$ d'\^etre invariant sous le groupe d'automorphisme ext\'erieur de $<H>$ provenant de $G$. Ensuite Arthur, pour $G$ connexe 
quasi-d\'eploy\'e, montre que l'application:
$$
f\in I_{cusp}(G) \mapsto \oplus_{<H>}f^{<H>}
$$
induit une isom\'etrie de $I_{cusp}(G)$ sur $\oplus_{<H>}I_{cusp}^{st, OUT_{G}(<H>)}(H)$ le produit scalaire sur $I_{cusp}(G)$ est le produit scalaire ordinaire tandis que le produite scalaire sur le membre de droite est la somme pond\'er\'ee des produits scalaires ordinaires, les coefficients \'etant l'analogue local des $i(G,H)$ de Langlands-Kottwitz. La d\'emonstration de \cite{arthurselecta} s'appuie sur \cite{waldspurgertransfert} 1.2. Le point cl\'e est de d\'emontrer la surjectivit\'e:
$$
I_{cusp}(G) \rightarrow I_{cusp}^{st,OUT_{G}(<H>)}(H).
$$ 
Pour le cas qui nous int\'eresse ici, $OUT$ est trivial; cela facilite la d\'emonstration.
Ceci a \'et\'e repris dans le cadre non connexe par \cite{waldspurger1} V mais o\`u seul le groupe endoscopique principal (celui qui contr\^ole la stabilit\'e) est trait\'e. C'est cette d\'emonstration dont on va v\'erifier qu'elle se g\'en\'eralise. Comme l'a remarqu\'e Arthur le point de d\'epart est le transfert au niveau des alg\`ebres de Lie c'est \`a dire la formule
 \cite{waldspurgertransfert} 1.2, maintenant valable dans le cas de $\tilde{G}$ et $<H>$ une donn\'ee endoscopique g\'en\'erale gr\^ace \`a \cite{waldspurgernouveau}. 
 
Il y a 2 points dans cette d\'emonstration; le premier point consiste \`a d\'efinir $I_{cusp}^{<H>-st}(\tilde{G})$ comme sous-espace de $I_{cusp}(\tilde{G})$; ce sont l'ensemble des fonctions sur les int\'egrales orbitales elliptiques qui  \`a l'int\'erieur d'une classe de conjugaison stable se transforment via le facteur de transfert relatif \`a la donn\'ee $<H>$. Et il faut d\'emontrer la d\'ecomposition (en tenant compte du fait que $OUT$ est trivial):
$$
I_{cusp}(\tilde{G})=\oplus_{<H>}I_{cusp}^{<H>-st}(\tilde{G}).
$$
Le deuxi\`eme point consiste \`a montrer que le transfert identifie $I_{cusp}^{<H>-st}(\tilde{G})$ et $I_{cusp}^{st}(H)$. Pour ce deuxi\`eme point on peut reprendre la d\'emonstration formelle d'Arthur \cite{arthurselecta} puisque $OUT$ est trivial ou celle de \cite{waldspurger1} VI.1 en y rempla\c{c}ant $stable$ par $<H>-stable$.
 
Pour le premier point, on fixe $<H>$ une donn\'ee endoscopique elliptique de $\tilde{G}$ et $I$ un 
\'el\'ement $<H>-$stable de $I_{cusp}(\tilde{G})$, on doit montrer qu'il existe $f\in C_{cusp}(\tilde{G})$ tel que $I$ soit la fonction associ\'ee aux int\'egrales orbitales de $f$.

Un argument de partition de l'unit\'e ram\`ene \`a faire cette d\'emonstration localement pr\`es des points semi-simples elliptiques. 

Fixons un \'el\'ement semi-simple elliptique de $\tilde{G}$ ou plus exactement sa classe de conjugaison stable. On a donc un ensemble fini de classes de conjugaison \`a l'int\'erieur de cette classe de conjugaison stable et il faut contruire la fonction cherch\'ee au voisinage de chacun de ces points de 
fa\c{c}on compatible \`a la conjugaison stable. On est donc amen\'e \`a travailler au voisinage de 0 dans l'alg\`ebre de Lie du centralisateur de chacun de ses points et la premi\`ere difficult\'e est que les centralisateurs ne sont pas isomorphes mais sont des formes int\'erieures l'un de l'autre. La remarque qui permet de travailler a \'et\'e faite en \cite{waldspurger1} V.3; le torseur qui ram\`ene \`a la forme quasid\'eploy\'ee de ces centralisateurs est compatible \`a la conjugaison stable (cf. le d\'ebut de V.3 dans loc. cite). Cette remarque est reprise dans le cadre g\'en\'eral en \cite{waldspurgernouveau} paragraphe 3.

On est donc ramen\'e \`a travailler avec la forme quasid\'eploy\'ee de ces centralisateurs mais il faut encore suivre les facteurs de transfert; si la classe stable de vient pas de $<H>$, on prend \'evidemment la fonction $0$ et il n'y a rien \`a faire. Sinon on est ramen\'e \`a la ''d\'efinition d'une donn\'ee endoscopique'' \cite{waldspurgernouveau} 3.5; c'est la g\'en\'eralisation au cas non connexe de l'utilisation de Kottwitz-Shelstad (\cite{ks}) par Arthur dans \cite{arthurselecta} et s'appuie fortement sur les d\'efinitions de la conjugaison stable de Labesse (\cite{labesseJIMJ}). Et on a l'\'egalit\'e des facteurs de transfert n\'ecessaire en \cite{waldspurgernouveau} 3.9. 
 
A c\^ot\'e de ce qui pr\'ec\`ede la d\'emonstration  au voisinage de l'origine pour des groupes assez g\'en\'eraux est compl\`etement limpide; c'est \cite{waldspurger1} V.2 o\`u il faut remplacer $stable$ par $<H>-stable$.


\subsection{Identit\'e de caract\`ere\label{caractere}}
Le r\'esultat principal de \cite{arthurselecta} consiste \`a montrer que pour qu'une identit\'e de 
caract\`eres obtenus sur les points elliptiques pour des repr\'esentations elliptiques se prolonge en une identit\'e de caract\`ere. Ceci reste vrai au moins pour ''stable'' dans la situation pr\'esente d'apr\`es 
\cite{waldspurger1}VI.3.

\bf Proposition. \sl Soit $\tilde{\Pi}$ une repr\'esentation virtuelle combinaison lin\'eaire de repr\'esentations  elliptiques de $G^+_{n}$. Soit $<H>$ une donn\'ee endoscopique elliptique de $\tilde{G}$ et $\Pi_{<H>}$ une repr\'esentation virtuelle de $<H>$ combinaison lin\'eaire de repr\'esentations elliptiques. Supposons que $tr\, \tilde{\Pi}$ soit un transfert de $tr\, \Pi_{<H>}$ sur les \'el\'ements elliptiques 
r\'eguliers. Alors $tr\, \tilde{\Pi}$ est un transfert de $tr\, \Pi_{<H>}$.\rm

\

On reprend \cite{waldspurger1}V.5: on globalise les groupes, pour cela il n'y a pas de probl\`eme; en suivant \cite{waldspurger1}, on note $u$ la place qui nous int\'eresse. On fixe deux places, $v_{1},v_{2}$ du corps global diff\'erentes de $u$ et en ces places on ne regardera que des fonctions cuspidales. On fixe encore une troisi\`eme place finie, $v_{3}$ diff\'erente des places d\'ej\`a utilis\'ees; en cette place on fixe une fonction cuspidale $<H>$-stable (cf. la preuve de \ref{decomposition}). On peut  utiliser la formule des traces simples d'Arthur gr\^ace \`a un r\'esultat de Mezo (\cite{mezo}) comme expliqu\'e en loc.cite.

On est alors dans une situation \'etudi\'ee en particulier par Labesse et on sait stabiliser la partie elliptique du c\^ot\'e g\'eom\'etrique de la formule des traces gr\^ace \`a \cite{labesseJIMJ}.
Avec les fonctions test satisfaisant aux 3 propri\'et\'es ci-dessus, seule la donn\'ee endoscopique $<H>$ intervient. Soient donc $f_{\tilde{G}}$ et $f_{H}$  des fonctions test sur $\tilde{G}$ et $H$ respectivement cuspidales aux 2 places fix\'ees ci-dessus et \`a la troisi\`eme place $f_{\tilde{G}}$ cuspidale et $<H>$-stable tandis que $f_{H}$ est stable;  avec la formule des traces simplifi\'ees on montre alors que $f_{H}$ est un transfert de $f_{\tilde{G}}$ si et seulement si la trace de $f_{\tilde{G}}$  est \'egale \`a la trace pour $f_{H}$; c'est la d\'emonstration faite dans loc. cite. En suivant toujours cette 
r\'ef\'erence (\cite{waldspurger1} IV.5 (3)) on sait aussi projeter les composantes \`a la place $u$ de 
$f_{\tilde{G}}$ et $f_{H}$ dans l'ensemble des fonctions cuspidales; ceci est d\'efini par le fait que sur toute repr\'esentation elliptique la trace de la fonction et celle de sa projection co\"{\i}ncident. La 
d\'emonstration consiste ensuite \`a montrer (en utilisant 2 fois l'\'equivalence ci-dessus) que si $f_{H}$ est un transfert de $f_{\tilde{G}}$ alors ceci reste vrai en rempla\c{c}ant les composantes \`a la place $u$ de ces fonctions test par leur projection cuspidale. En particulier les projections d\'efinies sont des fonctions cuspidales, celle sur $\tilde{G}$ \'etant un transfert de celle sur $H$. On note $f_{\tilde{G}}^u$ et $f_{H}^u$ les composantes en $u$ de $f_{\tilde{G}}$ et $f_{H}$ et $f_{cusp, \tilde{G}}^u$ et $f_{cusp,H}^u$ les projections. On a donc puisque $\pi$ et $\pi^H$ sont des combinaisons lin\'eaires de repr\'esentations elliptiques:
$$
tr\, \pi(f_{\tilde{G}}^u)=tr\, \pi(f_{cusp,\tilde{G}}^u); \qquad tr\, \pi^H(f_{H}^u)=tr\, \pi^H(f_{cusp,H}^u).
$$
On vient de voir que $f_{cusp,\tilde{G}}^u$ est un transfert de $f_{cusp,H}^u$.
Et l'\'egalit\'e $tr\, \pi(f_{\tilde{G}}^u)=tr\, \pi^H(f_{H}^u)$ r\'esulte donc de l'hypoth\`ese.
\subsection{Transfert et module de Jacquet\label{moduledejacquet}}
Soit $\tilde{\pi}$ une repr\'esentation irr\'eductible elliptique de $\tilde{G}_{n}$. En couplant les sections pr\'ec\'edentes, on trouve pour toute donn\'ee endoscopique elliptique $<H>$ de $\tilde{G}_{n}$, une repr\'esentation virtuelle $\pi^{<H>}$ combinaison lin\'eaire de repr\'esentations elliptiques de $H$ 
telle que l'on ait pour tout \'el\'ement  $\gamma\in \tilde{G}_{n,reg}$
$$
tr\, \tilde{\pi}(\gamma)=\sum_{<H>}\sum_{\gamma_{H}\in H} \Delta^{<H>}(\gamma_{H},\gamma) \, tr\, \pi^{<H>}(\gamma_{H}).\eqno(1)
$$Dans la somme de droite les \'el\'ements de $\gamma_{H}$ parcourt un ensemble de repr\'esentants de classes de conjugaison. A priori cette \'egalit\'e n'est vraie que pour les \'el\'ements elliptiques mais comme expliqu\'e ci-dessus, elle s'\'etend \`a tout \'el\'ement fortement r\'egulier.

\

On fixe $M$ un sous-groupe de Levi de $GL(n,E)$ de la forme $GL(a,E)\times GL(n-2a,E)\times GL(a,E)$ et un parabolique de Levi M; on prend le parabolique triangulaire par blocs sup\'erieur. On consid\`ere un point $\gamma=z \gamma' $ o\`u $\gamma'\in \tilde{M}$ et o\`u $z$ est un \'el\'ement de $E^*$ identifi\'e au centre du premier $GL(a,E)$. On applique (1) \`a un tel $\gamma$ pour $z$ dilatant les sous-groupes du radical unipotent de $P$, $\theta$-stables. On note $res_{P}\tilde{\pi}$ la restriction de $\pi$ au parabolique $P$ mais que l'on voit comme une repr\'esentation de $\tilde{M}$ puisque $\theta$ continue d'agir. Pour $z$ suffisamment dilatant, on obtient en suivant Casselman une \'egalit\'e:
$$
tr\, \tilde{\pi}(z\gamma')=\delta^{1/2}(z\gamma') tr\, res_{P}\tilde{\pi}(z\gamma'),
$$
o\`u $\delta$ est la fonction module. Il faut faire un calcul analogue pour le membre de droite; on fixe $<H>$ une donn\'ee endoscopique. On \'ecrit $H=U(n_{1},E/F) \times U(n_{2},E/F)$. On se limite aux \'el\'ements $\gamma'$ qui sont dans $M$ et dont la composante sur $GL(a,E)\times 1 \times GL(a,E)$ est elliptique. On note $M_{1}^H$ et $M^H_{2}$ les sous-groupes de Levi de $H$ isomorphes respectivement \`a $GL(a,E)\times U(n_{1}-2a,E/F)\times U(n_{2},E/F)$ et $U(n_{1},E/F)\times GL(a,E)\times U(n_{2}-2a,E/F)$ et $P_{i}^H$ les sous-groupes paraboliques ''standard'' correspondant; $M_{1}^H$ n'existe que si $n_{1}\geq 2a$ et de m\^eme pour $M_{2}^H$. Soit $\gamma''_{H}$ un \'el\'ement de $H$ tel que $\Delta^{<H>}(\gamma''_{H},z\gamma')\neq 0$. Alors $\gamma''_{H}$ est conjugu\'e d'un \'el\'ement de $M_{1}^H$ ou d'un \'el\'ement de  $M_{2}^H$, les 2 n'\'etant pas exclusifs;   on fixe $i=1,2$ tel que $\gamma''_{H}$ soit conjugu\'e d'un \'el\'ement de ce Levi (on consid\`ere alors $\gamma''_{H}$ comme un \'el\'ement de ce Levi) et on \'ecrit $\gamma''=z\gamma'_{H}$ o\`u $z$ est vu comme un 
\'el\'ement du facteur $GL(a,E)$ du Levi dans lequel $\gamma''_{H}$ est elliptique. On calcule encore 
$tr\, \pi^{<H>}(z\gamma'_{H})$ \`a l'aide de la restriction de $\pi^H$ au parabolique du Levi fix\'e.
On obtient donc un calcul de $tr\, res_{P}\tilde{\pi}(z\gamma')$ en fonction de ces restrictions. Les deux termes d\'ependent de $z$ dans un c\^one; cette \'egalit\'e se prolonge \`a tout $z$ et on peut donc faire $z=1$. Fixons encore $<H>$ comme ci-dessus et $\gamma'_{H}$, et notons $i=1,2$ l'indice tel que $\gamma'_{H}$ soit elliptique dans $M_{i}^H$. Il est clair que $M_{i}^H$ est naturellement une donn\'ee endoscopique elliptique pour $M$ et le point est de comparer les facteurs de transfert. C'est une situation simple puisque que l'on passe d'un groupe \`a l'un de ses Levi. Pr\'ecisons les notations; on \'ecrit $\gamma'=m h m'\theta$ avec $m,m'\in GL(a,E)$ et $h\in GL(n-2a,E)$ et on suppose que $\gamma'_{H}$ est  dans $M_{2}^H$; on \'ecrit alors $\gamma'_{H}=m_{2}h'_{H}$ o\`u $m_{2}\in GL(a,E)$ et $m_{2}\in U(n_{1},E/F)\times U(n_{2}-2a,E/F)$. On v\'erifie que si les classes stables de $\gamma'_{H}$ et $\gamma'$ se correspondent, alors il en est de m\^eme de $h\theta$ et $h'_{H}$ et $m_{2}$ est dans la m\^eme classe de conjugaison que $m\theta(m')$. Aux fonctions modules pr\`es qui disparaissent dans les calculs, on trouve, pour $\omega$ un caract\`ere convenable de $E^*$ (cf. introduction)
$$
\Delta^{<H>}(\gamma_{H},\gamma)=\omega(det(m_{2}))\Delta^{<U(n_{1}\times U(n_{2}-2a)>}(h'_{H},h\theta).
$$
Si on travaille avec le parabolique $P_{1}^H$, il n'y a pas le caract\`ere $\omega$. Le caract\`ere $\omega$ est un caract\`ere de $E^*$ qui prolonge le caract\`ere de $F^*$ correspondant via le corps de classe \`a l'extension $E$ de $F$.

 Finalement on trouve, avec les notations ci-dessus en y faisant $m'=1$ ce qui est loisible et avec $m$ elliptique dans $GL(a,E)$:
$$
tr \,  res_{P}\tilde{\pi}(mh\theta)=\sum_{<H>, i=1,2} \sum_{\gamma'_{H}}\Delta^{<M_{i}^H>}(h'_{H},h\theta) tr\,  res_{P_{i}^H} \pi^{<H>}(m h'_{H}) \omega^{i-1}(m).\eqno(2)
$$
On peut en tirer des r\'esultats plus pr\'ecis en d\'ecoupant suivant les repr\'esentations cuspidales du groupe $GL(a,E)$ qui interviennent. Il y  a des regroupements possibles et donc des simplifications possibles. 

 On consid\`ere 
$res_{P)}\tilde{\pi}$ comme une repr\'esentation de $GL(a,E)\times GL(n-2a,E)\times GL(a,E)$.
Et on d\'ecompose cette repr\'esentation dans le groupe de Grothendieck convenable sous la forme:
$$
\oplus_{\sigma,\sigma'}\sigma\otimes \pi(\sigma,\sigma')\otimes \sigma',
$$
o\`u $\sigma,\sigma'$ parcourt un ensemble de repr\'esentants des classes d'isomorphisme de 
repr\'esentations irr\'educ\-tibles de $GL(a,E)$ et o\`u $\pi(\sigma,\sigma')$ est une repr\'esentation virtuelle de $GL(n-2a,E)$. Supposons que $\sigma'\simeq \, ^\theta(\sigma)$, alors $\pi(\sigma,\, ^\theta(\sigma))$ a naturellement une action de $\theta$ et est donc une repr\'esentation virtuelle de $\tilde{G}_{n-2a}$, 
repr\'esentation que l'on \'ecrit $\tilde{\pi}(\sigma,^\theta(\sigma))$. De l'alg\`ebre lin\'eaire simple montre que si l'on \'ecrit :
$$
tr\, res_{P}\tilde{\pi}(mh\theta)=\sum_{\sigma}tr \, \sigma(m) tr\, \tilde{\pi}(\sigma,\theta(\sigma))(h\theta),$$
o\`u $\sigma$ parcourt un ensemble de repr\'esentant des classes d'\'equivalence de repr\'esentation 
irr\'educ\-tibles de $GL(a,E)$. De m\^eme pour toute repr\'esentation virtuelle $\Pi^{<H>}$, on \'ecrit, pour $i=1,2$:
$$
res_{P_{i}^H}\Pi^H=\sum_{\sigma} \sigma\otimes \Pi^{H}(\sigma,i).
$$
On veut d\'ecouper l'\'egalit\'e (2) suivant les repr\'esentations $\sigma$; comme on a suppos\'e $m$ elliptique dans $GL(a,E)$, c'est un peu d\'elicat. On se limite donc aux repr\'esentations $\sigma$ cuspidales (non n\'ecessairement unitaires). On obtient donc pour toute repr\'esentation cuspidale de $GL(a,E)$, $\sigma$, une \'egalit\'e de traces:
$$
tr\, \tilde{\pi}(\sigma,\theta(\sigma))(h\theta)=$$
$$
\sum_{<H=U(n_{1})\times U(n_{2})>}\sum_{h'_{H}}\Delta^{<U(n_{1}-2a)\times U(n_{2})>}(h'_{H},h\theta) tr\, \Pi^{H}(\sigma,1)
+
$$
$$
\sum_{<H=U(n_{1})\times U(n_{2})>}\sum_{h'_{H}}\Delta^{<U(n_{1})\times U(n_{2}-2a)>}(h'_{H},h\theta) tr\, \Pi^H(\omega\otimes \sigma,2). \eqno(3)
$$
On g\'en\'eralise ce r\'esultat en rempla\c{c}ant $\tilde{\pi}$ par une somme avec coefficients de 
repr\'esentations.

\

\bf Notations commode pour les modules de Jacquet.\rm

Il sera commode d'\'ecrire $\tilde{\pi}(\sigma,\theta(\sigma))(h\theta)=:Jac^\theta_{\sigma}\tilde{\pi}$. Pour $\tau$ une repr\'esentation de $U(n,E/F)$ et pour
 $\sigma$ une repr\'esentation cuspidale comme ci-dessus, on \'ecrira $Jac_{\sigma}\tau$, ce qui est \'ecrit ci-dessus $\tau(\sigma,1)$

\section{Appartenance \`a un paquet stable\label{appartenance}}
On fait d'abord une remarque g\'en\'erale. Soit $\tau'$ une s\'erie discr\`ete d'un groupe $U(m',E/F)$.
On note $f_{\tau'}$ l'\'el\'ement de $I_{cusp}(U(m',E/F)$ qui correspond \`a $\tau'$;  la projection de 
$I_{cusp}(U(m',E/F))$ sur $I^{st}_{cusp}(U(m',E/F))$ (par exemple suivant \cite{arthurselecta} 3.5) est non nulle;
l'argument  m'a \'et\'e donn\'e par Waldspurger. On regarde les germes du caract\`ere au voisinage de l'origine; ils se d\'eveloppent par degr\'e d'homog\'en\'eit\'e et le degr\'e formel est l'un de ces termes. Ce terme est stable et sa projection stable est donc non nulle et il y a donc un germe de la projection de $f_{\tau'}$ qui est non nul. D'o\`u la non nullit\'e de la projection de $f_{\tau'}$ sur $I_{cusp}^{st}(U(n,E/F)$
\section{Point de r\'eductibilit\'e des induites de cuspidales pour les groupes unitaires}
\subsection{Int\'egralit\'e\label{integralite}}
 Soit $\tau$ une repr\'esentation cuspidale de $U(n,E/F)$ et soient $\rho$ une repr\'esentation cuspidale irr\'educ\-tible $\theta$-invariante de $GL(d_{\rho},E)$. Soit $x\in {\mathbb R}_{>0}$ tel que  l'induite $\rho\vert\,\vert^x\times \tau$ soit r\'eductible.
 
 \
 
 \bf Proposition. \sl Avec les notations pr\'ec\'edentes, le r\'eel $x$ est un demi-entier ou encore $2x\in {\mathbb N}$.
\rm

\

Comme l'induite $\rho\vert\,\vert^x\times \tau$ est r\'eductible, il existe $\tau'$ une sous-repr\'esentation de cette induite qui est une s\'erie discr\`ete. Le module de Jacquet cuspidal de $\tau'$ est r\'eduit \`a l'unique terme $\rho\vert\,\vert^x\otimes \tau'$. On note $\Pi'$ un paquet stable de repr\'esentations contenant $\tau'$ et $\Sigma$ la repr\'esentation virtuelle image de $\Pi'$ dans la d\'ecomposition de 
$I_{cusp}(\tilde{G}_{n+2d_{\rho}})$. On calcule les modules de Jacquet des 2 membres et on projette sur la repr\'esentation cuspidale $\rho\vert\,\vert^x$ de $GL(d_{\rho},E)$ en suivant 
\ref{moduledejacquet} (3). Le seul groupe endoscopique intervenant ici est le $U(n+2d_{\rho})$-stable par construction mais ici $\Sigma$ est une repr\'esentation virtuelle. Quand on calcule les modules de Jacquet de  $\Pi'$,  le terme $\rho\vert\,\vert^x\otimes \tau$ ne peut dispara\^{\i}tre; en effet par 
r\'eciprocit\'e de Frobenius un tel terme ne peut provenir que de $\tau'$; ainsi il existe une 
repr\'esentation $\sigma'\in \Sigma'$ telle que
$$
Jac^\theta_{\rho\vert\,\vert^x}(\sigma')=\sigma'[\rho\vert\,\vert^x,\rho\vert\,\vert^{-x}]\neq 0.\eqno(1)
$$
On conna\^{\i}t la forme de $\sigma'$ c'est une induite de Steinberg de la forme
$$
\sigma'=\times_{(\rho',a')\in {\cal E}'}St(\rho',a').
$$
La non nullit\'e de (1) entra\^{\i}ne qu'il existe $(\rho',a')\in {\cal E}'$ tel que $\rho'=\rho$ et $x=(a'-1)/2$. D'o\`u la proposition.

\section{Finitude\label{finitude}}
Soit $\tau$ une repr\'esentation cuspidale de $U(n,E/F)$; on note $Red(\tau)$ l'ensemble des couples 
$(\rho,x_{\rho,\tau})$ tels que $\rho$ soit une repr\'esentation cuspidale irr\'eductible d'un groupe 
$GL(d_{\rho},E)$, $\theta$ invariante et $x_{\rho,\tau}$ soit un demi-entier strictement sup\'erieur \`a 1/2 tels que l'induite:
$$
\rho\vert\,\vert^{x_{\rho,\tau}}\times \tau
$$
soit r\'eductible. D'apr\`es \cite{silberger}, $x_{\rho,\tau}$ quand $\rho$ et $\tau$ sont fix\'es est uniquement d\'etermin\'e.

\

\bf Proposition. \sl L'ensemble $Red(\tau)$ est fini et on a l'in\'egalit\'e
$$
\sum_{(\rho,x_{\rho,\tau}\in Red(\tau)}(2x_{\rho,\tau}-1)d_{\rho}\leq n.
$$
\rm
Ceci est analogue \`a un r\'esultat de \cite{algebra} d\'emontr\'e par d'autres m\'ethodes; on am\'eliorera l'in\'egalit\'e dans la suite (de fa\c{c}on conforme \`a \cite{algebra}). Il suffit de montrer que l'in\'egalit\'e de la proposition est vraie pour tout sous-ensemble fini de $Red(\tau)$. On fixe ${\cal X}$ un 
sous-ensemble fini de $Red(\tau)$ et on consid\`ere $\tau_{{\cal X}}$ l'unique sous-module irr\'eductible de l'induite:
$$
\times_{(\rho,x_{\rho,\tau})\in {\cal X}}\rho\vert\,\vert^{x_{\rho,\tau}}\times \tau.
$$
On peut ordonner ${\cal X}$ comme on veut ci-dessus sans changer le r\'esultat et $\tau_{{\cal X}}$ est caract\'eris\'e par le fait que son module de Jacquet cuspidal est r\'eduit \`a l'ensemble des termes:
$$
\otimes_{(\rho,x_{\rho,\tau})\in {\cal X}}\rho\vert\,\vert^{x_{\rho,\tau}}\otimes \tau,\eqno(1)
$$
o\`u cette fois on consid\`ere tous les ordres possibles sur ${\cal X}$. On met $\tau_{{\cal X}}$ dans un paquet stable de repr\'esentations elliptiques et on consid\`ere l'image dans $I_{cusp}(\tilde{G}_{n_{{\cal X}}}$ pour $$n_{{\cal X}}:=n+\sum_{(\rho,x_{\rho,\tau})\in {\cal X}}d_{\rho}.
$$
On calcule les modules de Jacquet successifs en utilisant les repr\'esentations cuspidales $\rho\vert\,\vert^{x_{\rho,\tau}}$ associ\'ees aux \'el\'ements de ${\cal X}$; pour faire cela, on fixe un ordre sur ${\cal X}$ d'o\`u un terme dans le module de Jacquet de $\tau_{{\cal X}}$ du type (1). On v\'erifie que les seuls repr\'esentations elliptiques de $U(n_{{\cal X}},E/F)$ qui contiennent ce terme comme sous-quotient de leur module de Jacquet sont 
pr\'ecis\'ement $\tau_{{\cal X}}$. Ainsi $\tau$ est obtenu dans le calcul des modules de Jacquet successifs. Et il existe une repr\'esentation elliptique de $\tilde{G}_{n_{{\cal X}}}$, 
$\tilde{\pi}_{{\cal X}}$ dont les modules de Jacquets successifs calcul\'es comme expliqu\'e en 
\ref{moduledejacquet} (3) sont non nuls. 
On \'ecrit:
$$
\tilde{\pi}_{{\cal X}}=:\times_{(\rho',a')\in {\cal E}'}St(\rho',a').
$$
La non nullit\'e et le fait que les cuspidales intervenant dans ${\cal X}$ sont toutes diff\'erentes, 
entra\^{\i}nent qu'il existe un sous-ensemble ${\cal E}''$ de ${\cal E}'$ tel que
$$
{\cal X}=\{(\rho', (a'-1)/2); (\rho',a')\in {\cal E}''\}.
$$
On a $n_{{\cal X}}=\sum_{(\rho',a')\in {\cal E}'}a'd_{\rho'}$ tout simplement parce que 
$\tilde{\pi}_{{\cal X}}$ est une repr\'esentation de $\tilde{G}_{n_{{\cal X}}}$. Mais le module de Jacquet calcul\'e donne une repr\'esentation de $\tilde{G}_{n}$, d'o\`u
$$
n=n_{{\cal X}}-\sum_{(\rho',a')\in {\cal E}''}2d_{\rho'}=\sum_{(\rho',a')\in {\cal E}''}2d_{\rho'}(a'-3)/2+\sum_{(\rho',a')\in {\cal E}'-{\cal E}''}a'd_{\rho'}.
$$
Pour $(\rho',a')\in {\cal E}''$, on a $2(a'-3)/2=(a'-1)-1=2x_{\rho,\tau}-1$. D'o\`u encore
$$
n=\sum_{(\rho,x_{\rho,\tau})\in {\cal X}}d_{\rho}(2x_{\rho,\tau}-1)+\sum_{(\rho',a')\in {\cal E}'}a'd_{\rho'}.
$$
Cela donne a fortiori l'in\'egalit\'e cherch\'ee et d\'emontre la proposition.

\subsection{Propri\'et\'es g\'en\'erales des s\'eries discr\`etes et repr\'esentations elliptiques\label{proprietesdesseriesdiscretes}}

Ici on fixe $n$ et $\tau$ une repr\'esentation temp\'er\'ee.
Par des m\'ethodes totalement standard, on montre qu'il existe une repr\'esentation cuspidale $\tau_{0}$ d'un groupe $U(m_{0},E/F)$ avec $m_{0}\leq m$ et un ensemble totalement ordonn\'e de segments \`a la Zelevinski
$(\rho',a',b')\in {\cal J}$, o\`u $\rho'$ est une repr\'esentation cuspidale unitaire de $GL(d_{\rho'},E)$, et les $a',b'$ sont des r\'eels tel que $a'-b'+1 \in {\mathbb N}$ avec une inclusion:
$$
\tau \hookrightarrow \times_{(\rho',a',b')\in {\cal J}}<\rho\vert\,\vert^{a'}_{E},\rho\vert\,\vert^{b'}_{E}>\times \tau_{0}.\eqno(2)
$$
Ici je prends comme convention que $<\rho\vert\,\vert^{a'}_{E},\rho\vert\,\vert^{b'}_{E}>$ est la 
repr\'esentation $St(\rho',a'-b'+1)\vert\,\vert^{(a'+b')/2}$. Pour \'ecrire cela, on n'utilise rien sur $\tau$ sauf que c'est une repr\'esentation irr\'eductible. Rappelons quand m\^eme rapidement la d\'emonstration; on fixe d'abord un ensemble ${\cal J}'$, totalement ordonn\'e, form\'e de couples $\rho',y'$ tel que $\rho'$ soit une repr\'esentation cuspidale unitaire et $y'$ un r\'eel et tel que l'on ait une inclusion:
$$
\tau \hookrightarrow \times_{(\rho',y')\in {\cal J}'}\rho'\vert\,\vert^{y'}_{E}\times \tau_{0}.\eqno(3)
$$
On fixe un tel choix.  Ensuite on fixe $(\rho_{0},y_{0})$ tel que $y_{0}$ soit minimum et on pousse $\rho_{0}\vert\,\vert^{y_{0}}_{E}$ vers la gauche;  par exemple sans changer l'ensemble ${\cal J}'$ on change son ordre tel que (3) soit toujours satisfait mais $(\rho_{0},y_{0})$ est le plus ''petit'' possible (le plus \`a gauche possible). Cela permet de remplacer l'induite form\'ee par les repr\'esentations qui le 
pr\'ec\`edent et lui m\^eme par un segment comme ci-dessus.
On obtient alors la propri\'et\'e suivante suppl\'ementaire sur ${\cal J}$

(4) $\qquad$ soit $(\rho'',a'',b'')< (\rho',a',b')$ avec $\rho'\simeq \rho''$ pour l'ordre de ${\cal J}$, alors $b''\leq b'$. 

On peut obtenir une propri\'et\'e de plus quand on part d'un ensemble ${\cal J}'$ tel que parmi tous les choix possibles,  le nombre de $y'<0$ soit maximal pour le choix fait.
Cette propri\'et\'e de ${\cal J}'$ entra\^{\i}ne:  soit $(\rho',a',b')\in {\cal J}$ tel que $b'>0$, alors il existe 
$(\rho',a'',b'')\in {\cal J}$ (\'eventuellement le m\^eme \'el\'ement) avec $a''\geq a'$ et $\rho'\vert\,\vert^{b''}\times \tau_{0}$ est r\'eductible.

Si on utilise en plus l'hypoth\`ese que $\tau$ est une repr\'esentation temp\'er\'ee, on va v\'erifier que si 
${\cal J}$
v\'erifie (2) et (4) alors  pour tout $(\rho',a',b')\in {\cal J}$ n\'ecessairement $a'+b'\geq 0$.

En effet supposons  qu'il existe $(\rho,a',b')\in {\cal J}$ tel que $a'+b'< 0$; on en fixe le triplet de 
${\cal J}$ tel que $a'+b'$ soit minimal.  Sans changer ${\cal J}$ on fixe un ordre sur ${\cal J}$ tel que (2) soit encore vrai et le cas particulier de (4) quand on fixe notre $(\rho,a',b')$ et $(\rho',a',b')$ soit le plus ''petit'' possible pour cet ordre, il est plus simple de dire le plus \`a gauche possible dans (2). On va 
v\'erifier que $(\rho',a',b')$ est n\'ecessairement le premier \'el\'ement de ${\cal J}$. En effet soit  $(\rho'',a'',b'') \in {\cal J}$ pr\'ec\'edent $(\rho,a',b')$. Si 
la repr\'esentation induite:
$$
<\rho'',a'',b''>\times <\rho',a',b'>$$
est irr\'eductible, on peut \'echanger les 2 facteurs sans perdre l'inclusion en (2) et le cas particulier de (4) qui nous int\'eresse. Donc par minimalit\'e sur l'ordre de ${\cal J}$, cette induite est r\'eductible d'o\`u $\rho''\simeq \rho'$ et les segments $[a'',b'']$, $[a',b']$ sont li\'es; rappelons que ce sont des segments 
d\'ecroissants et que l'on a $b''\leq b'$ donc la liaison entra\^{\i}ne aussi $a''\in ]a',b'-1]$ et $b'' <b'$. 
D'o\`u $a''+b''<a'+b'$ ce qui contredit la minimalit\'e de $(\rho,a',b')$. On v\'erifie alors que l'exposant du module de Jacquet de $\tau$ donn\'e par r\'eciprocit\'e de Frobenius et (2) n'est pas dans le c\^one positif obtus ferm\'e comme il devrait l'\^etre (crit\`ere de Casselman) d'o\`u la contradiction cherch\'ee. Pour les s\'eries discr\`etes les in\'egalit\'es sont strictes.

\section{Changement de base r\'eciproque des repr\'esentations cuspidales}
Notre approche du changement de base se fait en utilisant les propri\'et\'es de r\'eductibilit\'e des induites de cuspidales; ces propri\'et\'es sont pr\'ecis\'ement control\'ees par des p\^oles de fonction $L$; c'est la th\'eorie d'Harish-Chandra amplifi\'ee par Shahidi.  
 \subsection{Support cuspidal \'etendu d'une s\'erie discr\`ete de $U(n,E/F)$ \label{supportcuspidal}}
 En \ref{integralite}, on  a montr\'e que les points de r\'eductibilit\'e des induites de cuspidales sont demi-entiers: pr\'ecis\'ement soit $\tau$ une repr\'esentation cuspidale de $U(n_{0},E/F)$ (ici $n_{0}$ peut \^etre $0$); soit $\rho$ une repr\'esentation cuspidale irr\'eductible d'un $GL(d_{\rho},E)$ et soit $x$ un 
 r\'eel strictement positif  tel que $\rho\vert\,\vert^{x}\times \tau_{0}$ soit r\'eductible. Alors il existe $a$ un entier naturel tel que $x=(a-1)/2$. Il faut remarquer que le cas $n_{0}=0$ est d\^u \`a Shahidi \cite{shahidi}.

 On a aussi montr\'e que pour un bon choix d'ensemble ${\cal J}$, on a l'inclusion 
 \ref{proprietesdesseriesdiscretes} (2)
 $$
 \tau \hookrightarrow \times _{(\rho',a',b')\in {\cal J}}<\rho,a',b'> \times \tau_{0}.
 $$On va v\'erifier encore que pour tout $(\rho',a',b')\in {\cal J}$, $a',b'$ sont des demi-entiers.

 Une telle assertion est \'evidemment ind\'ependante du choix de ${\cal J}$, c'est une propri\'et\'e de ''demi-int\'egralit\'e'' du support cuspidal de $\tau$. C'est \'el\'ementaire et on en rappelle la 
 d\'emonstration:

on fixe ${\cal J}$ un ensemble totalement ordonn\'e comme en  \ref{proprietesdesseriesdiscretes} 
v\'erifiant (2) et (4) de loc. cite. On impose en plus la propri\'et\'e de minimalit\'e
$$
\sum_{(\rho',a',b')\in {\cal J}}\sum_{x\in [a',b']}x \eqno(1)
$$
est minimal parmi les choix de ${\cal J}$ possible. Ceci est tout \`a fait compatible aux propri\'et\'es requises. Supposons qu'il existe $(\rho',a',b')\in {\cal J}$ avec $a',b'$ non demi-entiers. On peut supposer que cet \'el\'ement est le plus \`a droite possible. On peut le faire commuter avec tous les \'el\'ements plus \`a droite que lui car ceux l\`a font intervenir des demi-entiers et on a l'isomorphisme:
 $$
< \rho',a',b'>\times \tau_{0}\simeq <\rho',-b',-a'>\times \tau_{0}.
 $$Or on sait que $a'+b'\geq 0$ et $a'-b'+1\in {\mathbb N}$ (car $[a',b']$ est un segment d\'ecroissant). D'o\`u n\'ecessairement $a'+b'>0$ et  $\sum_{x\in [a',b']}x>\sum_{x'\in [-b',-a']}x'$. En ayant chang\'e $(a',b')$ en $(-b',-a')$, on va trouver un nouvel ensemble qui donnera un (1) plus petit que pour ${\cal J}$. Ceci est donc en contradiction avec la minimalit\'e fix\'ee.

\subsection{Pr\'eliminaire sur le support cuspidal \'etendu\label{preliminaire}}

Soit $\tau$ une repr\'esentation cuspidale de $U(n,E/F)$ et soit $\rho$ une repr\'esentation cuspidale 
irr\'eductible $\theta$ invariante de $GL(d_{\rho},E)$. On suppose qu'il existe $x>1/2$ tel que l'induite $\rho\vert\,\vert^x_{E}\times \tau$ soit r\'eductible; on sait d'apr\`es \cite{silberger} que $x$ est alors uniquement d\'etermin\'e et on pose $a_{\rho,\tau}:=2 x-1$. On sait maintenant que $a_{\rho,\tau}$ est un entier.

\

\bf Lemme. \sl Soit $\ell$ un demi-entier dans $  [1/2, (a_{\rho,\tau}+1)/2]$ tel que $2\ell+1$ soit de 
m\^eme parit\'e que $a_{\rho,\tau}$. Alors l'induite:
$$
\rho\vert\,\vert^{\ell}_{E}\times \cdots \times \rho\vert\,\vert^{(a_{\rho,\tau}+1)/2}_{E}\times \tau
$$
a un unique sous-module irr\'eductible et ce sous-module est une s\'erie discr\`ete.\rm

\

On note $m'$ l'entier $m+2(a_{\rho,\tau}+1)/2-\ell+1)d_{\rho}$. La repr\'esentation induite \'ecrite dans l'\'enonc\'e est une repr\'esentation de $U(m',E/F)$. On note $\sigma$ cette repr\'esentation induite.
L'unicit\'e du sous-module irr\'eductible r\'esulte de la r\'eciprocit\'e de Frobenius: soit $\tau_{0}$ un sous-module irr\'eductible. On a $Hom_{U(m',E/F)} (\tau_{0},\sigma)\neq 0$ entra\^{\i}ne que le module de Jacquet de $\tau_{0}$ admet le terme 
$$
\rho\vert\,\vert^{\ell}_{E}\otimes \cdots \otimes \rho\vert\,\vert^{(a_{\rho,\tau}+1)/2}\otimes \tau\eqno(1)
$$ 
comme quotient. Le terme que l'on vient d'\'ecrire intervient avec multiplicit\'e 1 comme sous-quotient du module de Jacquet de $\sigma$ et comme prendre les modules de Jacquet est un foncteur exact, 
$\tau_{0}$ est unique et intervient avec multiplicit\'e 1 comme sous-quotient de $\sigma$.

On pr\'ecise d'abord le lemme en annon\c{c}ant qu'en plus le module de Jacquet cuspidal de $\tau_{0}$ est r\'eduit au terme (1). On d\'emontre alors cette version plus forte du lemme par r\'ecurrence descendante sur $\ell$. Le r\'esultat est vrai pour $\ell=(a_{\rho,\tau}+1)/2$. On suppose donc $\ell<
(a_{\rho,\tau}+1)/2$; on change la notation $\tau_{0}$ pr\'ec\'edente en $\tau_{\ell}$ pour la faire 
d\'ependre de $\ell$. On consid\`ere \`a l'int\'erieur de $\sigma$ l'intersection des 2 sous-modules
$$
\rho\vert\,\vert^{\ell}_{E}\times \tau_{\ell+1}\cap <\rho,\ell, (a_{\rho,\tau}+1)/2>\times \tau,
$$
o\`u $<\rho,\ell, (a_{\rho,\tau}+1)/2>$ est l'unique sous-module irr\'eductible pour le groupe lin\'eaire convenable de l'induite:
$$
\rho\vert\,\vert^{\ell}_{E}\times \cdots \times \rho\vert\,\vert^{(a_{\rho,\tau}+1)/2}_{E}.
$$
Par unicit\'e de $\tau_{\ell}$ l'intersection contient $\tau_{\ell}$.
On calcule les modules de Jacquet de chacune des induites \'ecrites le module de Jacquet de $\tau_{\ell}$ est inclus dans l'intersection des modules de Jacquet. Il suffit donc de v\'erifier que cette intersection est r\'eduite au terme (1). On l'\'ecrit explicitement pour $\ell=(a_{\rho,\tau}-1)/2$, le cas 
g\'en\'eral suit la m\^eme ligne.
Le module de Jacquet de $\rho\vert\,\vert^{\ell}_{E}\times \tau_{\ell+1}$ a pour semi-simplifi\'e la somme des 4 termes (on pose $a:=a_{(\rho,\tau)}$) et on oublie l'indice $E$:
$$
\biggl(\rho\vert\,\vert^{(a-1)/2} \otimes \rho\vert\,\vert^{(a+1)/2}\otimes \tau\biggr) \oplus \biggl(
\rho\vert\,\vert^{(a+1)/2} \otimes \rho\vert\,\vert^{(a-1)/2}\otimes \tau \biggr) \oplus
$$
$$
\biggl(\rho\vert\,\vert^{-(a-1)/2} \otimes \rho\vert\,\vert^{(a+1)/2}\otimes \tau\biggr)\oplus\biggl(
\rho\vert\,\vert^{(a+1)/2} \otimes \rho\vert\,\vert^{-(a-1)/2}\otimes \tau\biggr).
$$
Le module de Jacquet de l'induite $<\rho,\ell, (a_{\rho,\tau}+1/1)>\times \tau$ est constitu\'e des 4 termes:
$$
\biggl(\rho\vert\,\vert^{(a-1)/2} \otimes \rho\vert\,\vert^{(a+1)/2}\otimes \tau\biggr) \oplus \biggl(
\rho\vert\,\vert^{(a-1)/2} \otimes \rho\vert\,\vert^{-(a+1)/2}\otimes \tau \biggr)\oplus
$$
$$
\biggl(\rho\vert\,\vert^{-(a+1)/2} \otimes \rho\vert\,\vert^{(a-1)/2}\otimes \tau\biggr)\oplus \biggl(
\rho\vert\,\vert^{-(a+1)/2} \otimes \rho\vert\,\vert^{-(a-1)/2}\otimes \tau\biggr).
$$
D'o\`u le r\'esultat dans ce cas particulier. Ce qui fait ''marcher'' la machine en g\'en\'eral est que le module de Jacquet de  la premi\`ere induite a dans tous ses termes $\rho\vert\,\vert^{(a+1)/2}$ 
tandis que la 2e,  a dans tous ses termes hormi le terme (1) $\rho\vert\,\vert^{-(a+1)/2}$. L'intersection est donc r\'eduite \`a (1).

\subsection{Support cuspidal \'etendu des repr\'esentations cuspidales des groupes unitaires\label{supportcuspidalcuspidal}}

\bf D\'efinition. \sl Soit $\tau$ une repr\'esentation cuspidale irr\'eductible de $U(m,E/F)$. Le support cuspidal \'etendu de $\tau$ est, par d\'efinition,  le support cuspidal pour le groupe lin\'eaire convenable de la 
repr\'esentation $$\pi_{\tau}:=\times_{\rho; a_{\rho,\tau}>0}\times_{b\in {\mathbb N}, b\leq a, b\equiv a[2]}St(\rho,b).$$
\rm

\

\bf Proposition. \sl (i)  $\pi_{\tau}$ est une repr\'esentation de $GL(n,E)$.

(ii)  Soit $\tilde{\pi}$ une repr\'esentation elliptique de $\tilde{G}_{n}$ dont la composante sur $I_{cusp}^{<U(n),st>}(U(n,E/F))$ contient $\tau$ de fa\c{c}on non nulle alors $\tilde{\pi}$ est un prolongement de $\pi_{\tau}$.\rm

\

On note $n'$ l'entier tel que $\pi_{\tau}$ soit une repr\'esentation de $GL(n',E)$. On montre d'abord que $n'\leq n$ en reprenant la d\'emonstration de \ref{finitude}.

Il faut g\'en\'eraliser la d\'efinition de $\tau_{\ell}$ de \ref{preliminaire} pour faire intervenir toutes les 
repr\'esentations $\rho$ telles que $a_{\rho,\tau}>0$; cet ensemble est fini gr\^ace \`a \ref{finitude}. On montre donc qu'il existe une s\'erie discr\`ete irr\'eductible $\tau'$ dont le module de Jacquet contient le terme:
$$\otimes_{\rho; a_{\rho,\tau}>0}
\biggl(
\rho\vert\,\vert^{\delta_{\rho}} \otimes \cdots \otimes \rho\vert\,\vert^{(a_{\rho,\tau}+1)/2}\biggr)\otimes \tau,
\eqno(1)$$o\`u $\delta_{\rho}$ vaut $1$ si $a_{\rho,\tau}$ est impair et $1/2$ sinon. Le module de Jacquet de $\tau'$ n'est pas irr\'eductible parce que l'on peut faire commuter des termes correspondant \`a des cuspidales diff\'erentes mais cela est la seule op\'eration possible. On note $m'$ l'entier tel que $\tau'$ soit une repr\'esentation de $U(m',E/F)$.

On va aussi avoir besoin de la propri\'et\'e d'unicit\'e suivante: soit $\tau''$ une repr\'esentation elliptique de $U(m',E/F)$ dont le module de Jacquet contient comme sous-quotient un des termes du type (1), alors $\tau''=\tau'$. En effet on applique la construction de \ref{proprietesdesseriesdiscretes} et on \'ecrit une inclusion:
$$
\tau'' \hookrightarrow \times_{(\rho',a',b')\in {\cal J}}<\rho',a',b'>\times \tau,
$$
o\`u les notations sont celles de loc.cite et $\tau$ est n\'ecessairement la repr\'esentation cuspidale 
fix\'ee. En comparant les supports cuspidaux ordinaires, on v\'erifie d\'ej\`a que pour tout $(\rho',a',b')\in {\cal J}$, $b'>0$; en effet il n'y a pas de multiplicit\'e dans le support cuspidal. Il faut maintenant v\'erifier les points suivants: soit $\rho'$ fix\'e et $(\rho',a',b')\in {\cal J}$. On suppose d'abord que $(\rho',a',b')$ ne pr\'ec\`ede aucun autre \'el\'ement de ${\cal J}$ de la forme $(\rho',a'',b'')$ alors $b'=(a_{\rho,\tau}+1)/2$.  S'il n'en \'etait pas ainsi on pourrait remplacer $b'$ par $-b'$ et contredire la positivit\'e d\'ej\`a 
d\'emontr\'ee. Ainsi pour cet \'el\'ement on a aussi $a'=b'$; puis on montre progressivement que si $(\rho',a',b')$ pr\'ec\`ede $t$ \'el\'ements de ${\cal J}$ de la forme $(\rho',a'',b'')$ alors $b'=(a_{\rho,\tau}+1)/2-t$: $b'$ ne peut \^etre plus grand car cela contredirait le fait que le support cuspidal ordinaire n'a pas de multiplicit\'e (avec un argument par r\'ecurrence sur $t$) ni plus petit car sinon on pourrait remplacer $b'$ par $-b'$ apr\`es avoir fait commuter et contredire la positivit\'e d\'ej\`a prouv\'ee; puis on obtient $a'=b'$ simplement comme cons\'equence de la connaissance du support cuspidal ordinaire (en termes imag\'es, tous les $\rho'\vert\,\vert^{x}$ avec $x>b'$ ont d\'ej\`a \'et\'e utilis\'es). Cela prouve l'unicit\'e cherch\'ee.

\

On consid\`ere la projection de $\tau'$ sur $I_{cusp}^{st}(U(m',E/F))$; cela donne une repr\'esentation virtuelle $\omega'$ qui contient de fa\c{c}on non nulle $\tau'$. Et on consid\`ere l'image de cet \'el\'ement stable dans $I_{cusp}(\tilde{G}_{m'})$. Cela fournit une combinaison lin\'eaire de repr\'esentations \`a laquelle on peut appliquer \ref{moduledejacquet} (3). On calcule successivement les modules de Jacquet pour les repr\'esentations cuspidales apparaissant dans (1) \`a gauche de $\tau$ dans l'ordre \'ecrit en commen\c{c}ant pas la gauche. Le calcul du module de Jacquet pour $\omega'$ est une combinaison lin\'eaire de termes mais le terme (1) ne peut pas dispara\^{\i}tre car il n'appartient qu'\`a $\tau'$ d'apr\`es l'unicit\'e ci-dessus. Ainsi dans le membre de gauche de (3) appliqu\'e successivement on obtient aussi un terme non nul. Soit donc $\Omega$ une repr\'esentation de $\tilde{G}_{m'}$  qui donne un terme non nul quand on prend ces modules de Jacquet. On sait pas hypoth\`ese que $\Omega$ est elliptique donc de la forme:
$$
\Omega=\times_{(\rho',b')\in {\cal E}}St(\rho',b')
$$
o\`u ${\cal E}$ est un ensemble convenable de couples. Le module de Jacquet est facile \`a calculer. Il est nul si ${\cal E}$ ne contient pas les couples $${\cal E}':=\cup_{\rho;a_{\rho,\tau}>0} \cup_{b\leq a; b\geq 0, b\equiv a [2]}(\rho,b+2).$$ Le module de Jacquet cherch\'e est alors
$$
\pi':=\times_{(\rho',b')\in {\cal E}'}St(\rho',b'-2) \times _{(\rho',b')\in {\cal E}-{\cal E}'}St(\rho',b').
$$
Or $\pi_{\tau}$ n'est autre que $\times_{(\rho',b')\in {\cal E}'}St(\rho',b'-2)$. D'o\`u $n'\leq n$.

Montrons que $n'\geq n$: on fixe $\tilde{\pi}$ et on lui associe un \'el\'ement de $I_{cusp}(\tilde{G}_{n})$. On suppose que dans la d\'ecomposition de \ref{decomposition} la composante de cet \'el\'ement fait intervenir $\tau$ de fa\c{c}on non nulle. C'est possible. On \'ecrit 
$$\pi:=\times_{(\rho,b)\in {\cal J}}St(\rho,b)$$ pour des choix convenables. On fixe $\rho$ tel qu'il existe $b$ avec $(\rho,b)\in {\cal J}$ et on fixe 
m\^eme $b_{0}$ maximal avec cette propri\'et\'e. On a alors $a_{\rho,\tau}=b_{0}$: en effet on 
consid\`ere la repr\'esentation:
$$
\pi(\psi'):=\times_{(\rho',b')\in {\cal J}-\{(\rho,b_{0})\}}St(\rho',b') \times St(\rho,b_{0}+2)
$$
un prolongement de cette repr\'esentation \`a $\tilde{G}_{n+2d_{\rho}}$. C'est une repr\'esentation $\theta$-discr\`ete (elliptique) par maximalit\'e de $b_{0}$. Puis on consid\`ere  le module de Jacquet $Jac^\theta_{\rho\vert\,\vert^{(b_{0}+1)/2}_{E}}$ qui est un prolongement de $\pi(\psi)$. Dans la projection sur la partie stable, ce module de Jacquet fait intervenir $\tau$ ce qui entra\^{\i}ne que dans sa d\'ecomposition $\pi(\psi')$ fait intervenir soit dans $I_{cusp}^{st}(U(n+2d_{\rho},E/F)$ un sous-quotient de l'induite $\rho\vert\,\vert^{(b_{0}+1)/2}_{E}\times \tau$ soit dans $I_{cusp}^{U(n)\times U(2d_{\rho}}(U(n,E/F)\times U(2d_{\rho},E/F)$ un sous-quotient de la repr\'esentation 
$$\tau \otimes \biggl(  ind\,\omega\otimes \rho\vert\,\vert^{(b_{0}+1)/2}\biggr).$$ 
On v\'erifie comme ci-dessus que dans les 2 cas le $Jac_{\rho\vert\,\vert^{-(b_{0}+1)/2}}$ (tordu par $\omega$ dans le 2e cas) de ce sous-quotient est 
n\'ecessairement 0 car $Jac^\theta_{\rho\vert\,\vert^{-(b_{0}+1)/2}}\pi(\psi')=0$. D'o\`u la r\'eductibilit\'e de l'induite $\rho\vert\,\vert^{(b_{0}+1)/2}\times \tau$ dans le premier cas et de l'induite de $\omega\otimes \rho\vert\,\vert^{(b_{0}+1)/2}_{E}$ dans le 2e cas. Le 2e cas est alors exclu car $(b_{0}+1)/2\geq 1$ et par les points de r\'eductibilit\'e pour une telle induite sont soit $0$ soit $1/2$. D'o\`u le fait que $b_{0}=a_{\rho,\tau}$ comme annonc\'e.
Il faut encore d\'emontrer que si $\rho$ est fix\'e et qu'il existe $b,b'$ des entiers avec $(\rho,b)$ et $(\rho,b')$ dans ${\cal J}$ alors $b$ et $b'$ ont m\^eme parit\'e. S'il n'en est pas ainsi on prend pour $b_{0}$ le plus grand entier tel que $(\rho,b_{0})\in {\cal J}$ et $b_{0}$ de parit\'e diff\'erente de $a_{\rho,\tau}$. On refait la d\'emonstration ci-dessus car $\pi(\psi')$ est encore une repr\'esentation elliptique de $\tilde{G}_{n+2d_{\rho}}$. Alors l'induite:
$$
\rho\vert\,\vert^{(b_{0}+1)/2}_{E}\times \tau
$$
serait encore r\'eductible ce qui contredirait l'unicit\'e des points de r\'eductibilit\'e dans ${\mathbb R}_{\geq 0}$ \cite{silberger}. D'o\`u la contradiction cherch\'ee.

 On a alors clairement
$$
\sum_{(\rho,b)\in {\cal J}}bd_{\rho}=n \leq \sum_{(\rho,b)\in {\cal J}; b\equiv a_{\rho,\tau}, b\leq a_{\rho,\tau}}bd_{\rho}=n'.
$$
En fait on d\'emontre beaucoup plus, puisque l'\'egalit\'e force que si $(\rho,b)\in {\cal J}$ avec $b>2$ alors $(\rho,b-2)\in  {\cal J}$. En d'autres termes ${\pi}=\pi_{\tau}$. On a donc d\'emontr\'e aussi (ii).
\subsection{Support cuspidal \'etendu pour les s\'eries discr\`etes d'un groupe unitaire\label{supportcuspidaletendu}}
On a d\'efini en \ref{supportcuspidalcuspidal} le support cuspidal \'etendu pour les repr\'esentations cuspidales des groupes unitaires. Soit maintenant $\tau$ une s\'erie discr\`ete irr\'eductible de $U(n,E/F)$. On revient \`a \ref{proprietesdesseriesdiscretes} (2), d'o\`u l'existence d'un ensemble de triplets $(\rho',a',b')$ tels que $\rho'$ soit une repr\'esentations cuspidales irr\'eductibles $\theta$ invariantes d'un groupe lin\'eaire et $a',b'$ un segment d\'ecroissant. On a montr\'e que $a',b'$ sont des demi-entiers; on a aussi la repr\'esentation cuspidale $\tau_{0}$ que l'on peut appeler comme en \cite{europe} et \cite{ams} le support cuspidal partiel et on pose:
$$
Supp cusp (\tau)=\cup_{(\rho',a',b')\in {\cal J}} \cup_{y\in [a',b']}\{\rho'\vert\,\vert^{y}_{E}, \rho'\vert\,\vert^{-y}_{E}\} \cup Supp cusp (\tau_{0}).
$$
Cette d\'efinition ne d\'epend que du support cuspidal ordinaire de $\tau$ et est donc ind\'ependante du choix de ${\cal J}$: en effet \'ecrivons $\tau$ comme sous-quotient d'une induite de la forme $\times_{\lambda\in {\cal S}}\lambda \times \tau_{0}$ o\`u ${\cal S}$ est un ensemble de repr\'esentations cuspidales irr\'eductibles de groupes lin\'eaires convenables et $\tau_{0}$ est une repr\'esentation cuspidale  irr\'eductible d'un groupe unitaire comme ci-dessus. Alors
$$
Supp cusp (\tau)=\cup_{\lambda\in {\cal S}}\{\lambda, \theta(\lambda)\} \cup Supp cusp(\tau_{0}).
$$
La d\'efinition s'\'etend sans probl\`eme aux repr\'esentations temp\'er\'ees puiqu'une repr\'esentation temp\'er\'ee irr\'eductible est un sous-module d'une induite de s\'eries discr\`etes. On peut aussi v\'erifier qu'une repr\'esen\-tation elliptique est combinaison lin\'eaire de repr\'esentation temp\'er\'ees ayant toutes m\^eme support cuspidal \'etendu puisque ces repr\'esentations sont toutes dans une m\^eme induite de s\'eries discr\`etes. Ceci permet donc de d\'efinir aussi le support cuspidal \'etendu d'une 
repr\'esentation elliptique. Le support cuspidal d'une repr\'esentation temp\'er\'ee fait intervenir des repr\'esentations cuspidales de la forme $\rho\vert\,\vert_{E}^x$ o\`u $\rho$ est unitaire non n\'ecessairement $\theta$-invariante; on v\'erifie ais\'ement que pour une repr\'esentation elliptique ceci ne peut pas se produire, $\rho$ est n\'ecessairement $\theta$-invariant.

\

\bf Lemme. \sl Il existe une unique repr\'esentation temp\'er\'ee de $GL(n,E)$, $\pi_{\tau}$ telle que le support cuspidal de $\pi_{\tau}$ soit pr\'ecis\'ement $Supp cusp (\tau)$.\rm

\

Avant de d\'emontrer ce lemme remarquons que cela donne une limitation forte sur les ensembles de repr\'esentations cuspidales qui peuvent \^etre un support cuspidal \'etendu d'une repr\'esentation 
temp\'er\'ee de $U(n,E/F)$. En effet, une repr\'esentation temp\'er\'ee de $GL(n,E)$ est n\'ecessairement de la forme
$$
\times_{(\rho,a)\in {\cal E}}St(\rho,a),
$$
avec des notations que l'on esp\`ere \'evidentes et le support cuspidale de cette repr\'esentation n'est autre que
$$
\cup_{(\rho,a)\in {\cal E}}\cup _{u\in [(a-1)/2,-(a-1)/2]}\rho\vert\,\vert^u_{E}.
$$

\

Il suffit de d\'emontrer le lemme pour $\tau$ une s\'erie discr\`ete. Soit donc $\tau$ une s\'erie discr\`ete.

Le support cuspidal \'etendu tel que d\'efini est un ensemble de couples $\rho,x$ (\'ecrit $\rho\vert\,\vert^x_{E}$) o\`u $\rho$ est une repr\'esentation cuspidale unitaire $\theta$-stable et $x$ un r\'eel. Pour chaque $\rho$  fix\'ee, l'ensemble des $x$ qui appara\^{\i}t est  un ensemble de demi-entiers stable par multiplication par $-1$.  Le lemme est \'equivalent \`a prouver que cet ensemble est une r\'eunion de segments centr\'es en l'origine et il est clair que si ceci est vrai la d\'ecomposition en segments est uniquement d\'etermin\'ee; une autre fa\c{c}on de dire les choses est que pour $GL(n,E)$ les 
repr\'esentations temp\'er\'ees sont uniquement d\'etermin\'ees par leur support cuspidal. Il n'y a donc que l'existence \`a d\'emontrer.

Fixons $(\rho',a',b')\in {\cal J}$ et supposons que $b'\leq 0$; alors 
$\cup_{y\in [a',b']}\{\rho'\vert\,\vert^{y}_{E},\rho'\vert\,\vert^{-y}_{E}\}$ est le support cuspidal de la 
repr\'esentation $St(\rho',2a'+1)\times St(\rho',-2b'+1)$. Dans l'\'ecriture \ref{proprietesdesseriesdiscretes} (2), on peut remplacer l'induite:
$$
\times_{(\rho',a',b')\in {\cal J}; b'>0}<\rho',a',b'>\times \tau_{0}
$$
qui sont les termes les plus \`a droite par un sous-quotient $\tau'$ convenable. Avec l'hypoht\`ese loisible que l'on a construit ${\cal J}$ en partant d'un ensemble ${\cal J}'$ avec un nombre maximum de $\rho''\vert\,\vert^{y''}_{E}$ avec $y''<0$, on sait que le module de Jacquet de $\tau'$ ne contient que des termes de la forme $\otimes_{\rho'',y''>0}\rho''\vert\,\vert^{y''}_{E}\times \tau_{0}$.
Ainsi $\tau'$ est une s\'erie discr\`ete et c'est ce que l'on a appel\'e en \cite{europe} une s\'erie discr\`ete fortement positive, on en a d\'ej\`a parl\'e ci-dessus. Il est tr\`es facile de classifier ces s\'eries discr\`etes quand on   conna\^{\i}t les 
diff\'erents termes $a_{\rho,\tau_{0}}$. Pour ces s\'eries discr\`etes, on d\'emontre que pour tout $\rho$ 
repr\'esentation cuspidale $\theta$-invariante irr\'eductible d'un $GL(d_{\rho},E)$ telle que $a_{\rho,\tau_{0}}\geq 0$, il existe un ensemble totalement ordonn\'e, \'eventuellement vide de segments 
d\'ecroissants, de cardinal $t_{\rho,\tau'}\geq 0$ inf\'erieur ou \'egal \`a $[(a_{\rho,\tau_{0}})/2]+1$, de la forme:
$$
\{ [x_{i}, (a_{\rho,\tau_{0}}+1)/2-i+1]; i\in [1,t_{\rho,\tau'}]\},
$$
o\`u $x_{1} < \cdots <x_{t_{\rho,\tau'}}$, avec une inclusion:
$$
\tau'\hookrightarrow \times_{\rho; a_{\rho,\tau_{0}}\geq 0}\times _{i\in [1,t_{\rho,\tau'}]}<\rho,x_{i},a_{\rho,\tau_{0}-i+1}>\times \tau_{0}.
$$
La d\'emonstration est sans difficult\'e et on renvoie \`a loc.cit. pour les d\'etails. On a alors:
$$
\pi_{\tau'}=\times_{\rho; a_{\rho,\tau_{0}}> 0}\biggl( \times_{i\in [1,t_{\rho,\tau'}}St(\rho,2x_{i}+1)
\times_{a<a_{\rho,\tau_{0}}-2t_{\rho,\tau'}} St(\rho,a)\biggr)$$
$$ \times_{\rho; a_{\rho,\tau_{0}}=0; t_{\rho,\tau'}=1}St(\rho, 2x_{1}).
$$
Il faut remarquer que si $a_{\rho,\tau_{0}}=0$ alors $t_{\rho,\tau'}$ vaut soit $0$ soit $1$.
Cela termine la preuve.

\

\bf Proposition. \sl Soit $\tau$ une repr\'esentation elliptique de $U(n,E/F)$ et soit $\tilde{\pi}$ une 
repr\'esentation elliptique irr\'eductible de $\tilde{G}_{n}$. On suppose que $\tilde{\pi}$ dans la 
d\'ecomposition \ref{decomposition} fournit un \'el\'ement de $I_{cusp}^{<U(n),st>}(U(n))$ non 
orthogonal \`a $\tau$. Alors $\tau$ est une s\'erie discr\`ete et $\tilde{\pi}$ est un prolongement de 
$\pi_{\tau}$. En particulier $\pi_{\tau}$ est $\theta$-discr\`ete.\rm

\

On \'ecrit $\tilde{\pi}$ comme un prolongement de la repr\'esentation:
$$
\pi\simeq \times_{(\rho,a)\in {\cal E}}St(\rho,a),
$$
o\`u ${\cal E}$ est un ensemble de couples form\'es d'une repr\'esentation cuspidale irr\'eductible 
$\theta$ invariante $\rho$ et d'un entier $a$; cet ensemble est sans multiplicit\'e. On consid\`ere la 
repr\'esentation 
$$
\pi^+:\simeq \times_{(\rho,a)\in {\cal E}}St(\rho,a+2).
$$C'est une repr\'esentation de $GL(n^+,E)$ o\`u $n^+=n+2\sum_{(\rho,a)\in {\cal E}}d_{\rho}$.
On ordonne ${\cal E}$ de telle sorte que si $(\rho',a')$ pr\'ec\`ede $(\rho'',a'')$ alors que $\rho'\simeq \rho''$ alors $a'<a''$. 
On calcule les modules de Jacquet successifs $$Jac^{\theta}_{\rho\vert\,\vert_{E}^{(a+1)/2}}\cdots  Jac^{\theta}_{\rho'\vert\,\vert_{E}^{(a'+1)/2}} \tilde{\pi}^+$$ o\`u l'ordre \'ecrit ci-dessus et l'ordre inverse de celui de ${\cal E}$; c'est-\`a-dire que l'on commence par le premier \'el\'ement pour aller vers le dernier. Chaque op\'eration donne une
repr\'esentation elliptique de la forme:
$$
\times_{(\rho,a)\in {\cal E}_{\leq}}St(\rho,a) \times_{(\rho,a)\in {\cal E}_{>}}St(\rho,a+2),
$$
o\`u ${\cal E}_{\leq}$ est le sous-ensemble de ${\cal E}$ utilis\'e pour prendre les modules de Jacquet et ${\cal E}_{>}$ celui qui n'a pas encore \'et\'e utilis\'e. A la derni\`ere \'etape, on trouve $\tilde{\pi}$. On 
v\'erifie alors qu'il existe $\tau_{+}$ intervenant dans $I_{cusp}^{<U(n^+),st>}(U(n^+,E/F)$ et tel que 
$\tau$ soit un sous-quotient de $Jac_{\rho\vert\,\vert_{E}^{(a+1)/2}; (\rho,a)\in {\cal E}}\tau^+$ o\`u l\`a aussi les $(\rho,a)$ sont pris dans l'ordre indiqu\'e ci-dessus. Le point ici est d'utiliser \ref{moduledejacquet} (3) et de montrer que les termes que l'on obtient sur $<U(n),st>$ ne peuvent provenir que de la donn\'ee endoscopique $<U(n^+),st>$ et non d'une donn\'ee de la forme $<U(n_{1})\times U(n_{2})>$ avec $n_{2}>0$. Ceci est juste une propri\'et\'e du support cuspidal ordinaire; en effet s'il n'en \'etait pas ainsi, il existerait une repr\'esentation elliptique de $U(n_{2},E/F)$, not\'ee $\tau'$, dont le support cuspidal ordinaire serait un sous-ensemble de $\{\rho\vert\,\vert^{(a+1)/2}_{E}; (\rho,a)\in {\cal E}\}$, en particulier ce support cuspidal n'est form\'e que de repr\'esentations de groupes lin\'eaires. Les exposants sont en valeur absolue, sup\'erieurs ou \'egaux \`a 1. Avec \ref{proprietesdesseriesdiscretes} on v\'erifie que $\tau'$ est une s\'erie discr\`ete et pas seulement une repr\'esentation elliptique. Son support cuspidal \'etendu est alors par d\'efinition inclus dans $\{\rho\vert\,\vert^{(a+1)/2},\rho\vert\,\vert^{-(a+1)/2}; (\rho,a)\in {\cal E}$; mais un tel support n'est pas une union de segments centr\'es en $0$ et ne peut donc \^etre le support cuspidal d'une repr\'esentation temp\'er\'ee d'un groupe lin\'eaire. Ceci contredirait le lemme ci-dessus.

On a donc prouv\'e l'existence de $\tau^+$. Par construction le support cuspidal ordinaire de $\tau^+$ se d\'eduit donc de celui de $\tau$ en ajoutant l'ensemble $\{\rho\vert\,\vert^{(a+1)/2}_{E}; (\rho,a)\in {\cal E}$ d'o\`u
$$
Supp cusp (\tau^+)=Supp cusp(\tau) \cup_{(\rho,a)\in {\cal E}}\{(\rho\vert\,\vert^{(a+1)/2},\rho\vert\,\vert^{-(a+1)/2}\}.
$$
Puisque $\tau^+$ est une repr\'esentation elliptique, on sait aussi que son support cuspidal \'etendu s'\'ecrit sous la forme:
$$
\cup_{(\rho',a')\in {\cal E}'}\cup_{y'\in [(a'-1)/2,-(a'-1)/2]}\rho'\vert\,\vert^{y'}_{E}.
$$
On revient \`a \ref{proprietesdesseriesdiscretes} pour savoir simplement que dans l'inclusion de cette 
r\'ef\'erence, les  segments de ${\cal J}$ (avec les notations de loc. cite) sont n\'ecessairement de la forme $<\rho',(a'-1)/2,b'>$ avec l'existence de $(\rho',a')\in {\cal E}'$ et soit $b'>0$ soit $b'$ est 
lui-m\^eme de la forme $-(a''-1)/2$ avec $(\rho',a'')\in {\cal E}$ (cf. la preuve du lemme ci-dessus).

Ainsi quand on applique les formules standard de Bernstein-Zelevinski pour calculer les modules de Jacquet, on a $Jac_{\rho'\vert\,\vert^x_{E}}\tau^+=0$ pour $x>0$ sauf s'il existe $(\rho',a')\in {\cal E}'$ avec $x=(a'-1)/2$. De plus pour un tel terme la filtration correspondante est r\'eduite \`a un \'el\'ement: il faut changer le segment qui est soit de la forme $<\rho',(a'-1)/2,b'>$ et $<\rho',(a'-1)/2-1,b'>$ soit de la forme $<\rho',A',-(a'-1)/2>$ en $<\rho',A',-(a'-1)/2+1>$. Ainsi ce terme n'intervient plus pour les autres modules de Jacquet \`a cause de l'ordre mis. Et on montre progressivement que la non nullit\'e du module de Jacquet de $\tau^+$ entra\^{\i}ne que
 pour chaque $(\rho,a)\in {\cal E}$, il existe $(\rho',a')\in {\cal E}'$ avec $\rho'\simeq \rho$ et 
$a'=a+2$. 

En d'autres termes pour tout $(\rho,a)\in {\cal E}$, il existe $(\rho,a+2)\in {\cal E}'$; on note ${\cal E}''$ l'ensemble restant. Quand on passe de $Supp cusp (\tau^+)$ \`a $Supp cusp( \tau)$ on enl\`eve 
pr\'ecis\'ement les extr\^emit\'es des segments non dans ${\cal E}''$ et on ne touche pas aux \'el\'ements de ${\cal E}''$; les couples non dans ${\cal E}''$ fournissent donc au support cuspidal \'etendu de $\tau$ exactement le support cuspidal de $\pi$ et il ne reste donc plus qu'\`a d\'emontrer que ${\cal E}''$ est vide. Mais, on a, puisque ${\cal E}$ donne le support cuspidal de $\pi$
$$
n= \sum_{(\rho,a)\in {\cal E}}ad_{\rho}.
$$
D'o\`u $n^+=\sum_{(\rho,a)\in {\cal E}}(a+2)d_{\rho}=\sum_{(\rho',a')\in {\cal E}'-{\cal E}''}a'd_{\rho'}.
$
Et par d\'efinition du support cuspidal \'etendu de $\tau^+$ 
$$
n^+=\sum_{(\rho',a')\in {\cal E}'}a'd_{\rho'}.
$$
D'o\`u ${\cal E}''=\emptyset$; ceci a plusieurs cons\'equences. On sait que le support cuspidal \'etendu de $\tau$ correspond \`a un ensemble de segments centr\'es \`a l'origine, ensemble sans multiplicit\'e. Ainsi $\tau$ est une s\'erie discr\`ete et pas seulement une repr\'esentation elliptique. La fin de la preuve est alors claire.

\subsection{Premi\`ere description des paquets stables de s\'eries discr\`etes pour les groupes unitaires\label{paquetstable1}}
Soit $\tau$ une s\'erie discr\`ete irr\'eductible de $U(n,E/F)$; on a d\'efini le support cuspidal \'etendu de $\tau$ que l'on peut voir comme le support cuspidal de la repr\'esentation $\pi_{\tau}$ d\'ej\`a d\'efinie. 
R\'eciproquement soit $\pi$ une repr\'esentation $\theta$-discr\`ete de $GL(n,E)$, on note $\Pi(\pi)$ l'ensemble des s\'eries discr\`etes de $U(n,E/F)$ ayant comme support cuspidal \'et\'endu le support cuspidal de $\pi$.

\

\bf Proposition. \sl Soit $\pi$ une s\'erie $\theta$-discr\`ete de $GL(n,E)$; on suppose que $\Pi(\pi)$ est non vide. Il existe une unique (\`a un scalaire pr\`es) combinaison lin\'eaire stable de s\'eries discr\`etes irr\'eductibles de $U(n,E/F)$ de la forme:
$$
\sum_{c_{\tau},\tau\in {\Pi(\pi)}}c_{\tau}\tau.
$$
De plus pour tout $\tau\in {\Pi(\pi)}$, $c_{\tau}\neq 0$.\rm

\

Soit $\tau\in \Pi(\pi)$; on sait que la projection, $f_{\tau}^{st}$ de $\tau$ sur $I_{cusp}^{st}(U(n,E/F))$ est non nulle et il existe donc une combinaison lin\'eaire d'\'el\'ements, $f$, de $I_{cusp}(\tilde{G}_{n})$ qui correspond \`a $f_{\tau}^{st}$ dans \ref{decomposition}. Il y a donc au moins une repr\'esentation elliptique $\tilde{\pi}$ intervenant dans $f$ telle que sa d\'ecomposition fasse intervenir $\tau$. On a vu que $\tilde{\pi}$ est un prolongement de $\pi_{\tau}=\pi$. Ainsi la projection de $\pi_{\tau}$ sur 
$I_{cusp}^{st}(\tilde{G}_{n})$ donne une distribution stable non nulle de $U(n,E/F)$. On sait que toutes les s\'eries discr\`etes de $U(n,E/F)$ intervenant dans cette d\'ecomposition ont comme support cuspidal \'etendu celui de $\pi_{\tau}$, 
c'est-\`a-dire sont dans $\Pi(\pi)$. D'o\`u l'existence d'une distribution stable ayant les propri\'et\'es de l'\'enonc\'e; cette distribution v\'erifie $c_{\tau}\neq 0$ pour le $\tau$ que nous avons fix\'e. Montrons maintenant l'unicit\'e annonc\'ee dans l'\'enonc\'e. Soit donc $f^{st},f^{',st}\in I_{cusp}^{st}(U(n,E/F)$ dont le support est inclus dans $\Pi(\pi)$. On note $f,f'$ les images de ces \'el\'ements dans 
$I_{cusp}(\tilde{G}_{n})$. En utilisant d'abord une repr\'esentation $\tau$ dans le support de $f^{st}$ puis une repr\'esentation $\tau'$ dans le support de $f^{'st}$, on v\'erifie comme ci-dessus qu'un prolongement de  $\pi$ intervient dans $f$ et dans $f'$ avec un coefficient non nul; changer de prolongement revient \`a changer les signes dans $f^{st}$ ou $f^{'st}$ et on peut donc supposer que c'est le m\^eme prolongement, $\tilde{\pi}$, qui intervient. Quitte \`a multiplier $f^{st}$ par un nombre complexe, on peut supposer que $f-f'$ n'a plus $\tilde{\pi}$ dans son support. L'image de $f-f'$ dans \ref{decomposition} est pr\'ecis\'ement $f^{st}-f^{'st}$. Si $f^{st}-f^{'st}$ \'etait non nul, il faudrait que $f-f'$ ait un coefficient non nul sur $\tilde{\pi}$ ce qui est exclu. D'o\`u $f^{st}=f^{'st}$ comme annon\c{c}\'e.

Comme pour chaque \'el\'ement de $\Pi(\pi)$ il existe $f^{st}$ satisfaisant les hypoth\`eses de l'\'enonc\'e avec $c_{\tau}\neq 0$, l'unicit\'e que l'on vient de d\'emontrer assure la propri\'et\'e de non nullit\'e, 
$c_{\tau}\neq 0$ pour tout $\tau\in \Pi(\pi)$. Cela termine la preuve.

\subsection{Support cuspidal et p\^ole de fonctions L\label{supportcuspidaletfonctionsL}}
Soit $\rho$ une repr\'esentation cuspidale irr\'eductible, $\theta$-invariante de $GL(d_{\rho},E)$. On note $L(\rho,r_{A}',s)$ la fonction $L$ d'Asai-Shahidi de $\rho$; cette fonction $L$ s'introduit naturellement 
gr\^ace aux travaux de Shahidi\cite{shahidi} compl\'et\'es par ceux de Goldberg \cite{goldberg}, c'est elle qui contr\^ole 
l'irr\'eductibilit\'e de l'induite de $\rho$ au groupe unitaire quasid\'eploy\'e $U(2d_{\rho},E/F)$ \`a partir du parabolique de Levi $GL(d_{\rho},E)$. Cette induite est irr\'eductible si et seulement si $L(\rho,r_{A}',s)$ a un p\^ole en $s=0$. Si cette induite est irr\'eductible alors l'induite dans la m\^eme situation de $\rho\vert\,\vert^{1/2}_{E}$ est, elle, r\'eductible. Et vice et versa.

\

\bf Proposition. \sl Soit $\tau$ une s\'erie discr\`ete irr\'eductible de $U(n,E/F)$ et $\rho$ comme ci-dessus; on suppose qu'il existe un entier $a$ tel que $\rho\vert\,\vert^{(a-1)/2}_{E}$ est dans le support cuspidal \'etendu de $\tau$. Alors $a$ est pair si $L(\rho,r_{A}',s)$ a un p\^ole en $s=0$ et impair sinon.\rm

\

Une autre fa\c{c}on de dire les choses peut-\^etre plus parlante est la suivante: on a d\'efini $\pi_{\tau}$ comme repr\'esentation $\theta$-discr\`ete de $GL(n,E)$ telle que son support cuspidal est exactement le support cuspidal \'etendu de $\tau$. On \'ecrit $$\pi_{\tau}=\times_{(\rho,a)\in {\cal E}}St(\rho,a).$$
Et la proposition dit que si $(\rho,a)\in {\cal E}$ alors $a$ est pair si et seulement si $L(\rho,r_{A}',s)$ a un p\^ole en $s=0$. 

\

On suppose d'abord que $\tau$ est cuspidal; on note $x_{0}$ l'unique r\'eel positif ou nul tel que 
$\rho\vert\,\vert^{x_{0}}\times \tau$ soit r\'eductible. On sait d\'ej\`a que $2x_{0}+1$ est un entier et on va montrer qu'il est impair si et seulement si $L(\rho,r_{A}',s)$ n'a pas  de p\^ole en $s=0$. Cela suffit 
gr\^ace \`a \cite{europe} appendice que l'on ne r\'ecrit pas ici. On fait remarquer au lecteur que l'on 
d\'emontre ici m\^eme un peu plus que ce qui est annonc\'e dans l'\'enonc\'e puisque l'on caract\'erise tous les points de r\'eductibilit\'e des induites de cuspidales.

\

Le  caract\`ere $\omega$ non trivial de $E^*$,  dont la restriction \`a $F^*$ est le caract\`ere de $F^*$ correspondant par la th\'eorie du corps de classe \`a l'extension $E$ de $F$ est celui qui s'introduit dans les choix pour l'endoscopie.

Supposons d'abord que $L(\rho,r_{A}',s)$ n'a pas de p\^ole en $s=0$. En suivant \cite{goldberg} (introduction, (1)) on sait que c'est alors 
$L(\omega\otimes \rho,r_{A}',s)$ qui a un p\^ole en $s=0$. Ainsi l'induite pour $U(2d_{\rho},E/F)$ de la 
repr\'esentation $\omega\otimes \rho\vert\,\vert^{1/2}_{E}$ contient une sous-repr\'esentation 
irr\'eductible, $\tau_{2}$ qui est une s\'erie discr\`ete. On note $f_{\tau}^{st}$ la projection de $\tau$ sur $I_{cusp}^{st}(U(n,E/F))$ et $f^{st}_{2}$ la projection de $\tau_{2}$ sur $I_{cusp}^{st}(U(2d_{\rho},E/F))$. On consid\`ere l'\'el\'ement $f^{st}_{\tau}\otimes f^{st}_{2}$ de 
$$I_{cusp}^{<U(n)\times U(2d_{\rho})>}(U(n,E/F)\times U(2d_{\rho},E/F))$$ 
qu'il d\'efinit et l'image r\'eciproque $\cal{F}$ dans $I_{cusp}(\tilde{G}_{n+2d_{\rho}})$ de cet \'el\'ement (cf \ref{decomposition}). On rappelle que l'on identifie un \'el\'ement de $I_{cusp}^{<H>}$ et une  
repr\'esentation virtuelle de $H$ et on utilise la notation pour les modules de Jacquet donn\'ee \`a la fin de \ref{moduledejacquet}

On applique \ref{moduledejacquet} (3) \`a ces \'el\'ements pour $\sigma=\rho\vert\,\vert^{1/2}_{E}$. On veut v\'erifier que le terme de droite de loc. cite (3) contient l'\'el\'ement $f^{st}_{\tau}$ comme composante sur $I_{cusp}^{<U(n),st>}(U(n,E/F))$; pour cela, il suffit de v\'erifier que $f^{st}_{2}$ contient $\tau_{2}$ et des repr\'esentations $\tau'$ telles que 
$Jac_{\omega\rho\vert\,\vert^{1/2}_{E}}\tau'=0$. Soit $\tau'$ une repr\'esentation irr\'eductible de 
$U(2d_{\rho},E/F)$ dont le module de Jacquet contient comme sous-quotient la repr\'esentation 
$\omega\rho\vert\,\vert^{1/2}_{E}$ du sous-groupe de Levi $GL(d_{\rho},E)$; un tel $\tau'$ est certainement un sous-quotient de l'induite de $\omega\rho\vert\,\vert^{1/2}_{E}$. Mais le seul sous-quotient de cette induite qui soit temp\'er\'ee est $\tau_{2}$ d'o\`u l'assertion cherch\'ee. Ainsi 
$$Jac^\theta_{\rho\vert\,\vert_{E}^{1/2}}{\cal{F}}$$ est non nul avec les notations de loc.cite. Il est imm\'ediat que $Jac^\theta_{\rho\vert\,\vert_{E}^{1/2}}{\cal{F}}$ est encore \`a support dans les
repr\'esentations $\theta$-discr\`etes (on enl\`eve un facteur $St(\rho,2)$ quand un tel facteur 
appara\^{\i}t et sinon on trouve 0). Et la projection de cet \'el\'ement sur $I_{cusp}^{st}(U(n,E/F))$ doit donc co\"{\i}ncider avec $f^{st}_{\tau}$. Cela veut dire que ${\cal F}$ a une composante non nulle sur 
$\pi_{\tau}\times St(\rho,2)$. On sait au d\'epart que $\pi_{\tau}\times St(\rho,2)$ est $\theta$-discr\`ete, cela veut dire que si l'on \'ecrit:
$$
\pi_{\tau}=\times_{(\rho',a')\in {\cal E}}St(\rho',a')
$$
l'ensemble ${\cal E}$ ne contient pas $(\rho,2)$. Or on sait que si $a_{\rho,\tau}\geq 1$ alors d'une part
$x_{\rho,\tau}=(a_{\rho,\tau}+1)/2$ et d'autre part si $a_{\rho,\tau}$ est pair, alors ${\cal E}$ contient 
$(\rho,2)$. On a donc montr\'e que si $a_{\rho,\tau}\geq 1$ c'est un entier impair et $2x_{\rho,\tau}+1$ est un entier pair. On veut aussi montrer que le cas o\`u $x_{\rho,\tau}=1/2$ est impossible sous nos 
hypoth\`eses. Il faut faire le m\^eme raisonnement mais en rempla\c{c}ant $\tau$ par l'unique sous-module irr\'eductible de l'induite $\tau':=\rho\vert\,\vert^{1/2}_{E}\times \tau$. La cuspidalit\'e de $\tau$ n'a servi \`a rien ci-dessus. Avec les hypoth\`eses, $\pi_{\tau'}$ est de la forme $St(\rho,2)\times \pi'$ avec 
$\pi'$ convenable et on obtient la contradiction cherch\'ee.

R\'eciproquement supposons que $(a_{\rho,\tau}+1)/2$ est entier non demi-entier, c'est-\`a-dire que soit $a_{\rho,\tau}=-1$ soit $a_{\rho,\tau}$ est un entier impair. On a d\'efini, $\pi_{\tau}$ et on consid\`ere 
$$\pi:= St(\rho,2)\times \pi_{\tau}.$$
C'est encore une repr\'esentation $\theta$ discr\`ete. On en fixe un prolongement $\tilde{\pi}$ \`a 
$G^+_{n+2d_{\rho}}$. On regarde la d\'ecomposition de l'image de $\tilde{\pi}$ dans 
$I_{cusp}(\tilde{G}_{n+2d_{\rho}})$ suivant  \ref{decomposition}. On note $f^{<H>}$ les diff\'erentes composantes quand $<H>$ parcourt l'ensemble des groupes endoscopiques elliptiques de 
$\tilde{G}_{n+2d_{\rho}}$. On calcule le module de Jacquet $Jac^\theta_{\rho\vert\,\vert^{1/2}}\pi$ et on obtient $\pi_{\tau}$. Cela permet de retrouver la d\'ecomposition de $\pi_{\tau}$ dans 
$I_{cusp}(\tilde{G})$ en utilisant \ref{moduledejacquet} (3); on ne s'interesse qu'\`a la composante suivant $I_{cusp}^{<U(n),st>}(U(n,E/F)$. Il faut donc consid\'erer $Jac_{\rho\vert\,\vert_{E}^{1/2}}f^{<U(n+2d_{\rho}),st>}$ (avec les notations de loc.cit.) et 
$Jac_{\omega\rho\vert\,\vert_{E}^{1/2}}f^{<U(n)\times U(2d_{\rho})>}$ o\`u le module de Jacquet porte sur les repr\'esentations de $U(2d_{\rho},E/F)$. Ce terme est facile \`a calculer:  il vaut 0 sauf si la 
repr\'esentation $\tau_{2}$ construite ci-dessus existe, c'est-\`a-dire si $L(\rho,r_{A}',s)$ n'a pas de p\^ole en $s=0$. Sous cette hypoth\`ese, on \'ecrit
$$
f^{<U(n)\times U(2d_{\rho})>}=f_{1}\otimes \tau_{2} \oplus f',
$$
o\`u $f'$ ne fait plus intervenir $\tau_{2}$ et 
$$Jac_{\omega\rho\vert\,\vert^{1/2}_{E}}f^{<U(n)\times U(2d_{\rho})>}=f_{1}.$$
Si $L(\rho,r_{A}',s)$ n'a pas de p\^ole en $s=0$, on a la r\'eciproque cherch\'ee. On raisonne donc par l'absurde en supposant que cette fonction $L$ a un p\^ole. Alors, n\'ecessairement
 $Jac_{\rho\vert\,\vert^{1/2}}f^{<U(n+2d_{\rho}),st>}$ est la projection de $\pi_{\tau}$ sur 
 $U^{st}(U(n,E/F))$. Cette projection contient $\tau$ et il existe donc une repr\'esentation dans le support de 
 $f^{<U(n+2d_{\rho}),st>}$ qui admet $\rho\vert\,\vert^{1/2}\otimes \tau$ comme sous-quotient dans son module de Jacquet. Cette repr\'esentation est une s\'erie discr\`ete et  l'induite $\rho\vert\,\vert^{1/2}\times \tau$ doit donc \^etre r\'eductible. Ceci est contradictoire avec le fait que $x_{\rho,\tau}$ est entier  non demi-entier. Cela termine la preuve.

\subsection{Paquets stables de s\'eries discr\`etes pour les groupes unitaires\label{paquetstable2}}
Soit $\psi$ un homorphisme de $W_{E}\times SL(2,{\mathbb C})$ dans $GL(n,{\mathbb C})$. On suppose que $\psi$ est semi-simple born\'e et continu au sens usuel. On suppose que la 
repr\'esentation $\psi$ est sans multiplicit\'e et qu'elle se d\'ecompose en somme de repr\'esentations 
irr\'eductibles, 
$$
\psi=\oplus_{(\rho,a)\in {\cal E}}\rho\otimes \sigma_{a}, \eqno(1)
$$
o\`u ici $\rho$ est une repr\'esentation irr\'eductible de dimension finie $d_{\rho}$ de $W_{E}$ et $a$ un entier, $\sigma_{a}$ \'etant l'unique repr\'esentation irr\'eductible de $SL(2,{\mathbb C})$ de dimension $a$. On dit que $\psi$ est $\theta$-discr\`ete si $\psi$ est sans multiplicit\'e et si toutes les 
repr\'esentations $\rho$ intervenant dans (1) sont invariantes sous l'action du compos\'e de $g\mapsto ^tg^{-1}$ et de la conjugaison dans $W_{E}$ venant de l'extension $E/F$.

On dit que $\psi$ est stable si dans sa d\'ecomposition (1), on a en plus que pour tout 
$(\rho,a)\in {\cal E}$, $a$ est pair si et seulement si $L(\rho,r_{A}',s)$ a un p\^ole en $s=0$; ici il est plus simple de dire que dans cette fonction $L$, $\rho$ est la repr\'esentation cuspidale de $GL(d_{\rho},E)$ correspondant \`a $\rho$ par la correspondance de Langlands. Bien s\^ur on peut interpr\'eter cette 
d\'efinition par le fait que $\psi$ se prolonge en un homomorphisme de $W_{F}\times SL(2,{\mathbb C})$ dans le groupe dual de $U(n,E/F)$. Mais on reste ici beaucoup plus \'el\'ementaire.

On rappelle encore que  la correspondance de Langlands associe \`a tout morphisme $\psi$ $\theta$-discret une repr\'esentation $\theta$-discr\`ete de $GL(n,E)$ que l'on note $\pi(\psi)$. On sait associer \`a $\pi(\psi)$ un paquet stable de s\'eries discr\`etes de $U(n,E/F)$ (qui \'eventuellement peut \^etre vide) en prenant la projection de l'image d'un prolongement de $\pi(\psi)$ dans $I_{cusp}(\tilde{G}_{n})$ sur 
$I_{cusp}^{<U(n),st>}(U(n,E/F))$. On note $\Pi(\psi)$ le paquet obtenu. On a montr\'e en 
\ref{paquetstable1} que $\Pi(\psi)$ est exactement form\'e des s\'eries discr\`etes de $U(n,E/F)$ dont le support cuspidal \'etendu co\"{\i}ncide avec le support cuspidal de $\pi(\psi)$. Il ne nous reste plus qu'\`a compl\'eter cette description en montrant que $\Pi(\psi)$ est non vide exactement quand $\psi$ est stable. C'est l'objet du th\'eor\`eme suivant:

\

\bf Th\'eor\`eme. \sl Il existe une bijection entre l'ensemble des paquets stables de s\'eries discr\`etes de $U(n,E/F)$ est l'ensemble des classes de conjugaison de morphismes $\theta$-discrets et stables de $W_{E}\times SL(2,{\mathbb C})$ dans $GL(n,{\mathbb C})$. La bijection est celle expliqu\'ee ci-dessus.\rm

\

Le th\'eor\`eme r\'esulte exactement de \ref{paquetstable1} et de \ref{supportcuspidaletfonctionsL}

\section{Changement de base}

Le but de cette partie est de d\'emontrer qu'\`a toute repr\'esentation $\theta$-discr\`ete $\pi$ de $GL(n,E)$ est associ\'ee une unique donn\'ee endoscopique elliptique $<H>$ de $\tilde{G}_{n}$ telle que l'image d'un prolongement de  $\pi$ dans $I_{cusp}(\tilde{G}_{n})$ soit exactement dans l'image de 
$I_{cusp}^{<H>}(H)$ dans la d\'ecomposition \ref{decomposition}. On va d\'emontrer cela, en montrant que $<H>$ est uniquement d\'etermin\'e par le support cuspidal de $\pi$.

Soit $n=n_{1}+n_{2}$ une d\'ecomposition de $n$, le couple $(n_{1},n_{2})$ est ordonn\'e. On en 
d\'eduit une donn\'ee endoscopique $<H>$ de $\tilde{G}_{n}$.

Soient $\tau_{1}$ une s\'erie discr\`ete de $U(n_{1},E/F)$ et $\tau_{2}$ une s\'erie discr\`ete de $U(n_{2},E/F)$; on les suppose irr\'eductible. On leur a associ\'e via le support cuspidal des repr\'esentations 
$\pi_{\tau_{1}}$ et $\pi_{\tau_{2}}$ de $GL(n_{1},E)$ et $GL(n_{2},E)$ respectivement. On pose:
$$
\pi_{\tau_{1}\otimes \tau_{2}}=\pi_{\tau_{1}} \times (\omega\otimes \pi_{\tau_{2}}),
$$
o\`u $\omega$ est comme ci-dessus un caract\`ere de $E^*$ dont la restriction \`a $F^*$ est le 
caract\`ere de $F^*$ correspondant \`a l'extension $E/F$.

\

\bf Lemme. \sl Soit $\pi$ une repr\'esentation $\theta$-discr\`ete de $GL(n,E)$ dont on fixe un prolongement $\tilde{\pi}$ \`a $\tilde{G}_{n}$. On note $f_{\pi}$ l'image de $\tilde{\pi}$ dans 
$I_{cusp}(\tilde{G}_{n})$. On suppose que la projection de $f_{\pi}$ sur $I_{cusp}^{<H>}(H)$ est non nulle. Alors $\pi=\pi_{\tau_{1}\otimes \tau_{2}}$.\rm

\

Ce lemme, dans le cas o\`u 
$n_{1}=n$, est un cas particulier de la proposition de  \ref{supportcuspidaletendu}; il s'en d\'eduit aussi si $n_{2}=n$, puisque la tensorisation par $\omega$ ram\`ene au cas stable. La d\'emonstration du cas g\'en\'eral suit le m\^eme principe. On fixe $\pi$ comme dans 
l'\'enonc\'e que l'on \'ecrit:
$$
\pi=\times_{(\rho,a)\in {\cal E}} St(\rho,a).
$$
On note $\pi^+:=\times_{(\rho,a)\in {\cal E}}St(\rho,a+2)$; c'est une repr\'esentation de $GL(n^+,E)$ o\`u $n^+$ est convenable. Et on ordonne encore ${\cal E}$ de fa\c{c}on \`a ce que si $(\rho',a')$ pr\'ec\`ede $(\rho'',a'')$ avec $\rho'=\rho''$ alors $a'<a''$. On fixe un prolongement $\tilde{\pi}^+$ de $\pi^+$ et on regarde l'image de $\tilde{\pi}^+$ dans $I_{cusp}(\tilde{G}_{n^+})$ et la d\'ecomposition de cet \'el\'ement suivant \ref{decomposition}. On va encore prendre les modules de Jacquet successifs suivant les 
\'el\'ements $\rho'\vert\,\vert^{(a'+1)/2}_{E}$ pour $(\rho',a')$ parcourant ${\cal E}$ dans l'ordre. A chaque \'etape on obtient une repr\'esentation elliptique et \`a la fin on trouve $\pi$. Quand on applique successivement \ref{moduledejacquet} (3), on obtient simplement l'existence d'une d\'ecomposition de 
${\cal E}$ en 2 sous-ensembles ${\cal E}_{1}\cap {\cal E}_{2}$, une donn\'ee endoscopique 
$<U(n_{1}^+)\times U(n_{2}^+)>$ de $\tilde{G}_{n^+}$ et un \'el\'ement $\tau_{1}^+\otimes \tau_{2}^+$ dans le support de la projection de $\tilde{\pi}^+$ sur 
$I_{cusp}^{<U(n_{1}^+)\times U(n_{2}^+)>}(U(n_{1}^+,E/F)\times U(n_{2}^+,E/F))$ tel que
$\tau_{1}^+$ ait dans son module de Jacquet un sous-quotient isomorphe \`a 
$\otimes_{(\rho',a')\in {\cal E}_{1}}\rho'\vert\,\vert^{(a'+1)/2}_{E}\otimes \tau_{1}$ et $\tau_{2}^+$ ait dans son module de Jacquet un sous-quotient isomorphe \`a $\otimes _{(\rho',a')\in {\cal E}_{2}}\omega\rho'\vert\,\vert^{(a'+1)/2}_{E}\otimes \tau_{2}$. 
On sait que $\tau_{1}$ est une s\'erie discr\`ete; en particulier si $(\rho',a')\in {\cal E}_{1}$ alors $a'$ est pair exactement si $L(\rho',r_{A}',s)$ a un p\^ole en $s=0$. De m\^eme si $(\rho',a')\in {\cal E}_{2}$ alors $a'$ est pair exactement si $L(\omega\rho',r_{A}',s)$ a un p\^ole en $s=0$; c'est la condition oppos\'ee \`a la pr\'ec\'edente puisqu'elle est \'equivalente \`a $L(\rho',r_{A}',s)$ n'a pas de p\^ole en $s=0$. Ainsi la d\'ecomposition de ${\cal E}$ en ${\cal E}_{1}\cup {\cal E}_{2}$ est compl\`etement d\'etermin\'ee a priori. Ensuite ce sont les arguments standard:
$$
Supp cusp (\tau_{1}^+)=\cup_{(\rho',a')\in {\cal E}_{1}}
\{\rho'\vert\,\vert^{(a'+1)/2},\rho'\vert\,\vert^{-(a'+1)/2}\} \cup Supp cusp (\tau_{1}).
$$
Cela prouve  encore que 
$$\pi_{\tau_{1}}=\times_{(\rho',a')\in {\cal E}_{1}}St(\rho',a') \times \pi'_{1},\eqno(1)$$ 
o\`u $\pi'_{1}$ est une repr\'esentation $\theta$-discr\`ete convenable. On fait la m\^eme chose en rempla\c{c}ant $1$ par $2$ et on trouve que $$\pi_{\tau_{2}}=\times_{(\rho',a')\in {\cal E}_{2}}\omega \otimes St(\rho',a') \times \pi'_{2}\eqno(2)$$ o\`u $\pi'_{2}$ est une repr\'esentation $\theta$-discr\`ete convenable. En comparant les dimensions, on trouve encore que $\pi'_{1}$ et $\pi'_{2}$ sont 
n\'ecessairement triviales. En mettant ensemble (1) et (2) tordu par $\omega$, on obtient
$$\pi_{\tau_{1}\otimes \tau_{2}}=
\pi_{\tau_{1}}\times (\omega\otimes\pi_{\tau_{2}})=\times_{(\rho',a')\in {\cal E}}St(\rho',a')=\pi.
$$
D'o\`u le lemme.

\subsection{Changement de base r\'eciproque\label{changementdebasereciproque}}
\bf Th\'eor\`eme. \sl Soit $\pi$ une repr\'esentation $\theta$ discr\`ete de $GL(n,E)$. Alors il existe une  donn\'ee endoscopique elliptique $<H>$ de $\tilde{G}_{n}$ tel que $\pi$ soit $<H>$-stable. Et $<H>$ est \'evidemment unique avec cette propri\'et\'e.\rm

\

Soit $\pi$ comme dans l'\'enonc\'e et $\tilde{\pi}$ un prolongement de $\pi$; on note $f_{\pi}$ l'image de $\pi$ dans $I_{cusp}(\tilde{G}_{n})$. Soit $<H>$ une donn\'ee endoscopique de $\tilde{G}_{n}$ telle que $f_{\pi}$ ait une projection non nulle sur $I_{cusp}^{<H>}(H)$. On \'ecrit $<H>=<U(n_{1}\times U(n_{2})>$ et on a vu qu'il existe $\tau_{1},\tau_{2}$ des s\'eries discr\`etes irr\'eductibles de $U(n_{1},E/F)$ et $U(n_{2},E/F)$ respectivement telles que $$\pi=\pi_{\tau_{1}}\times (\omega \otimes \pi_{\tau_{2}}).$$
On \'ecrit $\pi=\times_{(\rho',a')\in {\cal E}}St(\rho',a')$ et pour $i=1,2$, 
$$\pi_{\tau_{i}}=\times_{(\rho',a')\in {\cal E}_{i}}St(\rho',a').$$
On sait que $\pi=\pi_{\tau_{1}\otimes \tau_{2}}$ ce qui est \'equivalent \`a dire que
$$
{\cal E}={\cal E}_{1}\cup_{(\rho',a')\in {\cal E}_{2}}(\omega\rho',a').
$$
Or soit $(\rho',a')\in {\cal E}$; alors $(\rho',a')$ ne peut \^etre dans ${\cal E}_{1}$ que si $a'$ est pair quand $L(\rho',r_{A}',s)$ n'a pas de p\^ole en $s=0$ et impair sinon. Et $(\omega\rho',a')$ ne peut \^etre dans ${\cal E}_{2}$ que si $a'$ est pair quand $L(\omega\rho',r_{A}',s)$ a un p\^ole en $s=0$ et impair sinon. Les propri\'et\'es de la fonction $L(\rho',r_{A}',s)$ et la pariti\'e de $a'$ d\'etermine donc si $(\rho',a')\in {\cal E}_{1}$ ou si $(\omega\rho',a')\in {\cal E}_{2}$ les 2 possibilit\'es \'etant exclusives l'une de l'autre. D'o\`u l'unicit\'e de $<H>$ et le th\'eor\`eme. 

\section{Cardinal des paquets stables\label{cardinal}}
On suppose ici que $U(n,E/F)$ est quasid\'eploy\'e.

On fixe un morphisme $\psi$ de $W_{E}\times SL(2,{\mathbb C})$ dans $GL(n,{\mathbb C})$, $\theta$-discret; on a d\'efini le paquet de s\'eries discr\`etes $\Pi(\psi)$. On note $I_{cusp}^U[\psi]$ le sous-espace vectoriel de $I_{cusp}(U(n,E/F))$ engendr\'e par l'image des \'el\'ements de $\Pi(\psi)$. Soit $H$ un groupe endoscopique de $U(n,E/F)$; $H$ est de la forme $U(n_{1},E/F)\times U(n_{2},E/F)$ avec $n_{1}+n_{2}=n$ et l'ordre ici n'importe pas. En particulier si $n_{1}=n_{2}$, la donn\'ee endoscopique qui vient avec $H$ admet un automorphisme ext\'erieur provenant de $U(n)$. Pour $n_{1},n_{2}$ tel que $n_{1}+n_{2}=n$, on d\'efinit d'abord $I_{cusp}^{n_{1},n_{2}}[\psi]$ en g\'en\'eralisant la d\'efinition ci-dessus et en sommant sur tous les morphismes, $\psi'$ de $W_{E}\times SL(2,{\mathbb C})$ dans $GL(n_{1},{\mathbb C})\times GL(n_{2},{\mathbb C})$ dont l'image par l'inclusion de 
$GL(n_{1},{\mathbb C})\times GL(n_{2},{\mathbb C})$ dans $GL(n,{\mathbb C})$ est conjugu\'e de 
$\psi$.  Et on note $I_{cusp}^{n_{1},n_{2},st}[\psi]$ le sous-espace de $I_{cusp}^{n_{1},n_{2}}[\psi]$ 
form\'e des \'el\'ements stables.

Ici on ne peut plus ignorer le fait que certaines repr\'esentations elliptiques ne sont pas des s\'eries 
discr\`etes. On note $I_{cusp}^{n,ell}$ le sous-espace vectoriel de $I_{cusp}(U(n,E/F)$ engendr\'e par l'image des repr\'esen\-ta\-tions elliptiques qui ne sont pas des s\'eries discr\`etes. On a d\'ej\`a vu (cf. \ref{paquetstable2}) que $$I_{cusp}(U(n,E/F))=\oplus_{\psi}I_{cusp}[\psi]\oplus I_{cusp}^{n,ell},$$o\`u $\psi$ parcourt l'ensemble des morphismes $\theta$-discrets pris \`a conjugaison pr\`es somme de 
repr\'esentations stables.
$$
I_{cusp}^{st}(U(n,E/F))=\oplus _{\psi}{\mathbb C} f_{\psi},\eqno(1)$$
o\`u $f_{\psi}$ est la combinaison lin\'eaire stable de s\'eries discr\`etes dans $\Pi(\psi)$ (cf. encore 
\ref{paquetstable2}). On a en plus
$$
I_{cusp}^{n,ell}\cap I_{cusp}^{st}(U(n,E/F))=0,
$$puisque $I_{cusp}^{st}(U(n,E/F)=I_{cusp}^{st}(\tilde{G}_{n})$ a \'et\'e d\'ecrit par 
\ref{changementdebasereciproque} comme \'etant \'egal \`a (1).

Pour $n=n_{1}+n_{2}$ comme ci-dessus avec $n_{1}n_{2}\neq 0$, on d\'efinit diff\'eremment 
$I_{cusp}^{n_{1},n_{2},st,ell}$. On pose:
$$
I_{cusp}^{n_{1},n_{2},st,ell}:=\oplus_{\psi_{1}\times \psi_{2}}
I_{cusp}^{st}(U(n_{1},E/F)\times U(n_{2},E/F))[\psi_{1}\times \psi_{2}],
$$
o\`u la somme porte sur les morphismes $\psi_{i}$ (pour $i=1,2$) de $W_{E}\times SL(2,{\mathbb C})$ dans $GL(n_{i},{\mathbb C})$, $\theta$-discret mais tel que les repr\'esentations d\'efinies par $\psi_{1}$ et $\psi_{2}$ ne sont pas disjointes. En d'autres termes $\psi_{1}\times \psi_{2}$ vu comme 
repr\'esentation de dimension $n$ n'est pas $\theta$-discr\`ete.

 \
 
 On utilise la d\'ecomposition de \cite{arthurselecta} 3.5; c'est ici que l'on utilise le fait que $U(n,E/F)$ est quasid\'eploy\'e:
 $$
 I_{cusp}(U(n,E/F))=\oplus_{n_{1},n_{2}} I_{cusp}^{st, OUT}(U(n_{1},E/F)\times U(n_{2},E/F)), \eqno(1)
 $$
 o\`u la somme porte sur les couples non ordonn\'es $n_{1},n_{2}$ tel que $n=n_{1}+n_{2}$ et o\`u $OUT$ indique que si $n_{1}=n_{2}$ on consid\`ere les \'el\'ements invariants par l'action de l'automorphisme ext\'erieur  \'echangeant les 2 copies. Si $n_{1}\neq n_{2}$, $OUT$ est sans objet. Rappelons qu'ici on a suppos\'e $U(n,E/F)$ quasi-d\'eploy\'e.
 
 Soit $\psi$ comme ci-dessus et $n_{1}=n_{2}$;  l'espace $I_{cusp}^{n_{1},n_{2},st}[\psi]$ est stable sous-l'action de $OUT$. Son suppl\'ementaire naturel (l'espace des \'el\'ements instables, c'est-\`a-dire les int\'egrales orbitales dont la somme sur toute classe stable est nulle) est lui aussi stable sous $OUT$ et  on obtient facilement:
 $$
 I_{cusp}^{st, OUT}(U(n_{1},E/F)\times U(n_{2},E/F))=
 $$
 $$\sum_{\psi}I_{cusp}^{st, OUT}(U(n_{1},E/F)\times U(n_{2},E/F))[\psi] \oplus I_{cusp}^{st, ell, OUT}(U(n_{1},E/F)\times U(n_{2},E/F)).
 $$

 \
 
 \bf Lemme. \sl Soit $\psi$ un morphisme $\theta$-discret comme ci-dessus. Soit $n=n_{1}+n_{2}$. 
 
 L'image de $I_{cusp}^{st, OUT}(U(n_{1},E/F)\times U(n_{2},E/F))[\psi] $ dans $I_{cusp}(U(n,E/F))$ suivant (1) est incluse dans $I_{cusp}(U(n,E/F))[\psi]$.\rm
 
 \
 
On fixe $\tau\in \Pi(\psi)$ et on montre que la d\'ecomposition de $\tau$ vu comme \'el\'ement de 
$I_{cusp}(U(n,E/F))$ a une projection sur chaque $I_{cusp}^{st}(U(n_{1},E/F)\times U(n_{2},E/F))$ incluse dans $$I_{cusp}^{st}(U(n_{1},E/F)\times U(n_{2},E/F))[\psi].$$ 
On \'ecrit $f_{\tau}$ l'\'el\'ement de $I_{cusp}(U(n,E/F))$ correspondant \`a $\tau$. Et on d\'ecompose:
$$
f_{\tau}=\oplus_{H}f^H_{\tau},
$$
o\`u pour tout $H$ groupe endoscopique elliptique de $U(n,E/F)$, $f^H_{\tau}$ est un \'el\'ement de $I_{cusp}^{st}(U(n_{1},E/F)\times U(n_{2},E/F))$. On peut donc d\'ecomposer cet \'el\'ement suivant la base obtenue \`a l'aide du produit de morphisme $\psi_{1}\times \psi_{2}$, o\`u $\psi_{1}$ et $\psi_{2}$ sont 
$\theta$-discrets et on doit d\'emontrer que seuls interviennent les morphismes $\psi_{1}\times \psi_{2}$ qui sont conjugu\'es de $\psi$. C'est donc encore un probl\`eme de support cuspidal; il faut d\'emontrer que le support cuspidal \'etendu de $\tau$ est le support cuspidal de $\pi(\psi_{1})\times \pi(\psi_{2})$.

On reprend les m\'ethodes d\'ej\`a utilis\'ees; on \'ecrit 
$$
\pi(\psi)=\times_{(\rho,a)\in {\cal E}}St(\rho,a).
$$
On note $\psi^+$ le morphisme tel que 
$$
\pi(\psi^+)=\times_{(\rho,a)\in {\cal E}}St(\rho,a+2).
$$
Et on a montr\'e dans la preuve de la proposition de \ref{supportcuspidaletendu}, l'existence d'un \'el\'ement $\tau^+$ dans $\Pi(\psi^+)$  sous-module 
irr\'eductible de l'induite
$$
\times_{(\rho,a)\in {\cal E}} \rho\vert\,\vert^{(a+1)/2}_{E}\times \tau,
$$
o\`u ici l'ordre sur ${\cal E}$ est important; l'\'el\'ement $(\rho,a)$ pr\'ec\`ede (c'est -\`a-dire est \`a gauche ci-dessus) l'\'el\'ement $(\rho,a')$ si $a<a'$.

On \'ecrit la d\'ecomposition de $f_{\tau^+}$ comme ci-dessus et on calcule les modules de Jacquet par rapport aux \'el\'ements $\rho\vert\,\vert^{(a+1)/2}_{E}$ pris dans le m\^eme ordre: on commence par celui de gauche et on ''remonte''. Pour $H$ fix\'e, $\psi_{1}\times \psi_{2}$ et $\tau_{1} \otimes \tau_{2}$ une 
s\'erie discr\`ete de ce paquet, intervenant dans la d\'ecomposition de $f_{\tau}^H$, on trouve qu'il existe (au moins) un groupe endoscopique $H'_{1}\times H'_{2}$, (au moins)
un morphisme 
$\psi'_{1}\times \psi'_{2}$ et (au moins) une  s\'erie discr\`ete, $\tau'_{1,2}$ dans le paquet d\'efini, dont le module de Jacquet contient comme sous-quotient $$\biggl(\otimes_{(\rho',a')\in {\cal E}_{1}}\rho'\vert\,\vert^{(a'+1)/2}_{E}\otimes\tau_{1}\biggr)\otimes \biggl(\otimes _{(\rho'',a'')\in {\cal E}_{2}} \rho''\vert\,\vert^{(a''+1)/2}_{E}\otimes\tau_{2}\biggr),
$$o\`u la d\'ecomposition de ${\cal E}$ en ${\cal E}_{1}\cup {\cal E}_{2}$ d\'epend des choix et ne nous importe pas.
En prenant les supports cuspidaux \'etendus, on trouve:
$$
supp cusp (\pi(\psi'_{1})\times \pi(\psi'_{2})) = supp cusp (\pi(\psi_{1})\times \pi(\psi_{2}) \cup \{\rho\vert\,\vert^{(a+1)/2},\rho\vert\,\vert^{-(a+1)/2}\}.
$$
Ici on a \'et\'e un peu trop vite en faisant toutes les \'etapes en un seul coup mais elle sont progressives: ce que l'on veut est que les \'el\'ements $(\rho,(a+1)/2), (\rho,-(a+1)/2$ sont des extr\'emit\'es des segments de $\psi'_{1}\times \psi'_{2}$. On a d\'ej\`a vu dans la preuve de la proposition de \ref{supportcuspidaletendu} qu'il en \'etait bien ainsi. On a aussi vu qu'un simple calcul des dimensions assurent que tout segment dans le support cuspidal de $\pi(\psi'_{1})\times \pi(\psi'_{2})$ est obtenu ainsi. On en d\'eduit donc que le support cuspidal de $\pi(\psi'_{1})\times \pi(\psi'_{2})$ est pr\'ecis\'ement celui de la repr\'esentation:
$$
\times_{(\rho,a)\in {\cal E}}St(\rho,a+2).
$$
Et celui de $\pi(\psi_{1})\times \pi(\psi_{2})$ s'obtient en enlevant toutes les extr\'emit\'es et c'est donc celui de $\pi(\psi)$ comme annonc\'e.
\subsection{Cardinal\label{calculcardinal}}
\bf Th\'eor\`eme. \sl Soit $\psi$ un morphisme $\theta$-discret de $W_{E}\times SL(2,{\mathbb C})$ dans $GL(n,{\mathbb C})$. On note $\ell(\psi)$ la longueur de la repr\'esentation ainsi d\'efinie. Alors:
$$
\vert \Pi(\psi)\vert= 2^{\ell(\psi)-1}.
$$
\rm

\

On sait que le nombre d'\'el\'ements dans $\Pi(\psi)$ est exactement la dimension de l'espace vectoriel $I_{cusp}(U(n,E/F)[\psi]$. Et on sait aussi que $I_{cusp}^{st}[\psi]$ est de dimension 1. Soit $n_{1}+n_{2}=n$ une d\'ecomposition ordonn\'ee de $n$ et soit $\psi_{i}$ pour $i=1,2$ des morphismes de $W_{E}\times SL(2,{\mathbb C})$ dans $GL(n_{i},{\mathbb C})$. L'automorphisme de $I_{cusp}(U(n_{1},E/F)\times U(n_{2},E/F))[\psi_{1}\times \psi_{2}]$ sur $I_{cusp}(U(n_{2},E/F)\times U(n_{1},E/F))[\psi_{2}\times \psi_{1}]$ envoie la droite ''stable'' sur son homologue. On suppose que $\psi_{1}\times \psi_{2}$ est conjugu\'e de $\psi$; alors $\psi_{1}$ et $\psi_{2}$ n'ont aucune sous-repr\'esentation en commun. Cela prouve que si $n$ est pair, le groupe $OUT$ n'a pas de droite fixe dans $I_{cusp}^{st}(U(n/2,E/F)\times U(n/2,E/F))[\psi]$ d'o\`u
$$
dim\, I_{cusp}^{st,OUT}(U(n/2,E/F)\times U(n/2,E/F))=1/2 dim\, I_{cusp}^{st}(U(n/2,E/F)\times U(n/2,E/F))[\psi].
$$
On a donc
$$
dim\, I_{cusp}(U(n,E/F))[\psi]=1/2 \sum_{(n_{1},n_{2}); n_{1}+n_{2}=n}dim\, I_{cusp}^{st}(U(n_{1},E/F)\times U(n_{2},E/F))[\psi],\eqno(1)
$$
o\`u la somme porte sur les d\'ecompositions ordonn\'ees de $n$; entre autre la d\'ecomposition $(0,n)$ n'est pas la m\^eme que $(n,0)$. Soit $n_{1},n_{2}$ une d\'ecomposition de $n$. On a encore
$$
I_{cusp}^{st}(U(n_{1},E/F)\times U(n_{2},E/F))[\psi]=\oplus_{\psi_{1},\psi_{2}}I_{cusp}^{st}(U(n_{1},E/F))[\psi_{1}]\otimes I_{cusp}^{st}(U(n_{2},E/F))[\psi_{2}],
$$
o\`u $\psi_{1},\psi_{2}$ sont comme ci-dessus et $\psi_{1}\times \psi_{2}$ est conjugu\'e de $\psi$. Il faut donc calculer le nombre de telles d\'ecompositions. Pour faire ces calculs, on inverse les sommes: on 
d\'ecompose $\psi$ en le produit de 2 repr\'esentations de $W_{E}\times SL(2,{\mathbb C})$, produit ordonn\'e, et cette d\'ecomposition d\'etermine uniquement $n_{1}$ et $n_{2}$. Ces d\'ecompositions sont en bijection avec l'ensemble des applications de l'ensemble des sous-repr\'esentations 
irr\'eductibles incluses dans $\psi$ dans $\{\pm 1\}$. Il y a donc $2^{\ell(\psi)}$ telles d\'ecompositions. Il ne faut pas oublier le 1/2 qui vient de $OUT$ (cf (1) ci-dessus) et on trouve le r\'esultat.

\section{Classification des repr\'esentations cuspidales}
\subsection{D\'efinitions \label{definitions}}
On fixe un homomorphisme $\psi$ de $W_{E}\times SL(2,{\mathbb C})$ dans $GL(n,{\mathbb C})$, 
$\theta$ discret et on suppose que ce morphisme param\'etrise un paquet non vide, $\Pi(\psi)$ de 
s\'eries discr\`etes. Soit $s$ un \'el\'ement, $\theta$-invariant, du centralisateur de $\psi$ dans 
$GL(n,{\mathbb C})$. On note $z$ la matrice diagonale de $GL(n,{\mathbb C})$ de valeurs propres 
$-1$. Le centralisateur de $s$ dans $GL(n,{\mathbb C})$ est un produit de 2 groupes lin\'eaires 
$GL(n_{1},{\mathbb C}) \times GL(n_{2},{\mathbb C})$, o\`u $n_{1}$ est la dimension de l'espace propre correspondant \`a la valeur propre $-1$ de $s$. En rempla\c{c}ant $s$ par $zs$ on \'echange $n_{1}$ et $n_{2}$. Puisque $s$ est dans le centralisateur de $\psi$, le morphisme $\psi$ se factorise par 
$GL(n_{1},{\mathbb C})\times GL(n_{2},{\mathbb C})$ et on note $\psi_{s}$ cette factorisation. On obtient ainsi une repr\'esentation $\pi(\psi_{s})$ de $GL(n_{1},E)\times GL(n_{2},E)$. On la prolonge en une repr\'esentation de $\tilde{G}_{n_{1}}\times \tilde{G}_{n_{2}}$ ce qui donne naturellement un \'el\'ement de $I_{cusp}^{st}(\tilde{G}_{n_{1}}\times \tilde{G}_{n_{2}})$ et donc de $I_{cusp}^{st}(U(n_{1},E/F)\otimes I_{cusp}^{st}(U(n_{2},E/F)$. On note $\Psi_{s^.}$ la projection de cet \'el\'ement dans $I_{cusp}^{st,OUT}(H)$, o\`u $H$ est la donn\'ee endoscopique de $U(n,E/F)$ associ\'ee au groupe 
$U(n_{1},E/F)\times U(n_{2},E/F)$. Comme nous n'avons pas fix\'e le choix de l'extension, en travaillant avec $sz$ plut\^ot que $s$ on aurait trouv\'e le m\^eme r\'esultat au signe pr\`es. On fixe donc un choix et on note $f_{\psi,s}$ l'image de cet \'el\'ement dans $I_{cusp}(U(n,E/F))$. On a vu que cet \'el\'ement s'interpr\`ete comme une combinaison lin\'eaire des caract\`eres des repr\'esentations dans $\Pi(\psi)$. Ainsi par inversion, on 
d\'efinit pour tout \'el\'ement $\pi\in \Pi(\psi)$ des nombres complexes $d(s^.,\pi)$ uniquement 
d\'etermin\'es par l'\'egalit\'e de caract\`eres:
$$
tr\, \pi=\sum_{(s,sz)}d(s^.,\pi)f_{\psi,s}.
$$  
On d\'ecompose $\psi$ en repr\'esentations irr\'eductibles et soit $(\rho,a)$ un couple form\'e d'une 
repr\'esentation irr\'eductible de $W_{E}$ et d'un entier $a$ tel que, en notant $\sigma_{[a]}$ la 
repr\'esentation irr\'eductible de $SL(2,{\mathbb C})$ de dimension $a$, la repr\'esentation 
$\rho\otimes \sigma_{[a]}$ de $W_{E}\times SL(2,{\mathbb C})$ soit une sous-repr\'esentation de $\psi$. On suppose qu'il existe $0\leq b <a$ tel que, si $b\neq 0$,  $\rho\otimes \sigma_{[b]}$ soit aussi une sous-repr\'esentation de $\psi$ et si $b=0$ que $a$ est pair. On note alors $a_{-}$ le plus grand 
\'el\'ement $b$  v\'erifiant les 
propri\'et\'es ci-dessus. On pose $z_{\rho,a}$ l'\'el\'ement du centralisateur de $\psi$ dans 
$GL(n,{\mathbb C})$ dont les valeurs propres $-1$ ont exactement pour espace propre la somme de $\rho\otimes \sigma_{[a]}\oplus \rho\otimes \sigma_{[a_{-}]}$. On note $A(\psi)$ le sous-groupe du centralisateur de $\psi$ engendr\'e par ces \'el\'ements $z_{\rho,a}$ quand $(\rho,a)$ parcourt tous les couples possibles.

\bf Lemme. \sl Soient $(\rho,a)$, $z_{\rho,a}$ ayant les propri\'et\'es pr\'ec\'edentes. Pour tout $s$ comme ci-dessus, il existe un signe $\zeta_{s,\rho,a}$ tel que
$$
Jac_{\rho\vert\,\vert^(a-1)/2, \cdots, \rho\vert\,\vert^{(a_{-}+1)/2}} f_{\psi,s}=\zeta_{s,\rho,a}
Jac_{\rho\vert\,\vert^(a-1)/2, \cdots, \rho\vert\,\vert^{(a_{-}+1)/2}} f_{\psi,sz_{\rho,a}}.
$$
\rm
Il faut distinguer 2 cas. Dans le premier cas l'espace propre pour la valeur propre $-1$ de $s$ contient l'espace de la repr\'esentation $\rho\otimes \sigma_{[a]}$ et ne contient pas l'espace de la 
repr\'esentation $\rho\otimes \sigma_{[a_{-}]}$. Le deuxi\`eme cas est le cas o\`u la valeur propre $-1$ de $s$ a un espace propre qui contient la somme de ces 2 repr\'esentations. Quitte \`a changer $s$ en $sz$, on se trouve dans l'un ou l'autre cas.

Dans les 2 cas, la d\'emonstration utilise le fait que prendre les modules de Jacquet est compatible avec le transfert endoscopique. Pour simplifier l'\'ecriture, on enl\`eve les $_{E}$ des valeurs absolues. On a 
d\'efini $Jac^\theta_{\sigma}$ \`a la fin de \ref{moduledejacquet} pour $\sigma$ une repr\'esentation cuspidale; on g\'en\'eralise cette notation \`a un ensemble fini de repr\'esentations cuspidales, $\sigma_{1}, \cdots, \sigma_{k}$, en posant $Jac^\theta_{\sigma_{1}, \cdots, \sigma_{k}}:=Jac^\theta_{\sigma_{k}}\circ \cdots \circ Jac^\theta_{\sigma_{1}}$.

Consid\'erons d'abord le premier cas; on calcule $Jac^\theta_{\rho\vert\,\vert^{(a-1)/2}, \cdots
 \rho\vert\,\vert^{(a_{-}+1)/2}}\pi(\psi_{s})$; on \'ecrit $\psi_{s}=\psi_{1}\times \psi_{2}$ o\`u $\psi_{i}$ pour $i=1,2$ est \`a valeurs dans $GL(n_{i},{\mathbb C})$, d\'ecomposition suivant les espaces propres de $s$. On a alors facilement:
 $$
 Jac^\theta_{\rho\vert\,\vert^{(a-1)/2}, \cdots,
 \rho\vert\,\vert^{(a_{-}+1)/2}} \biggl(\pi(\psi_{1})\otimes \pi(\psi_{2})\biggr)
 $$
 $$
 = \biggl(Jac^\theta_{\rho\vert\,\vert^{(a-1)/2}, \cdots,
 \rho\vert\,\vert^{(a_{-}+1)/2}}\pi(\psi_{1})\biggr) \otimes \pi(\psi_{2}).
 $$
 On note $\psi'_{1}$ le morphisme qui se d\'eduit de $\psi$ en rempla\c{c}ant la sous-repr\'esentation 
 $\rho\otimes \sigma_{[a]}$ par $\rho\otimes \sigma_{[a_{-}]}$ et le r\'esultat est
 $\pi(\psi'_{1})\otimes \pi(\psi_{2})$. L'action de $\theta$ sur le r\'esultat est d\'etermin\'ee par le choix fait au d\'epart.
 
 On fait le m\^eme calcul en partant de $sz_{\rho,a}$. Alors le morphisme $\psi_{sz_{\rho,a}}$ vaut 
 $\psi'_{1}\otimes \psi'_{2}$ o\`u $\psi'_{1}$ est exactement comme ci-dessus et $\psi'_{2}$ se d\'eduit de 
 $\psi_{2}$ en rempla\c{c}ant la repr\'esentation $\rho\otimes \sigma_{[a_{-}]}$ en $\rho\otimes \sigma_{[a]}$. Et on obtient:
$$
 Jac^\theta_{\rho\vert\,\vert^{(a-1)/2}, \cdots,
 \rho\vert\,\vert^{(a_{-}+1)/2}} \biggl(\pi(\psi'_{1})\otimes \pi(\psi'_{2})\biggr)
$$
$$
 = \pi(\psi'_{1}) \otimes \biggl(Jac^\theta_{\rho\vert\,\vert^{(a-1)/2}, \cdots
 \rho\vert\,\vert^{(a_{-}+1)/2}} \pi(\psi'_{2})\biggr)
$$
$$
=\pi(\psi'_{1})\otimes \pi(\psi'_{2}).
$$Il n'y a pas de difficult\'e \`a faire agir le groupe des automorphismes ext\'erieurs venant de $U(n,E/F)$ sur la donn\'ee endoscopique pour garder l'\'egalit\'e au signe pr\`es des transferts endoscopiques $$Jac _{\rho\vert\,\vert^{(a-1)/2}, \cdots
 \rho\vert\,\vert^{(a_{-}+1)/2}} f_{\psi,s}$$ et $Jac_{\rho\vert\,\vert^{(a-1)/2}, \cdots
 \rho\vert\,\vert^{(a_{-}+1)/2}} f_{\psi,sz_{\rho,a}}$, le signe d\'epend de l'action de $\theta$;  ici c'est une \'egalit\'e dans $I_{cusp}(U(n-(a-a_{-})d_{\rho},E/F))$, on se trouve dans la partie elliptique non s\'erie discr\`ete.

 Consid\'erons maintenant le 2e cas: on note encore $\psi_{1}\times \psi_{2}$ la d\'ecomposition en produit de $\psi_{s}$. Ici on note $\psi'_{1}$ le morphisme qui se d\'eduit de $\psi_{1}$ en enlevant les 2 repr\'esentations $\rho\otimes \sigma_{[a]}$ et $\rho\otimes \sigma_{[a_{-}]}$. Et on a:
 $$
 Jac^\theta_{\rho\vert\,\vert^{(a-1)/2}, \cdots
 \rho\vert\,\vert^{(a_{-}+1)/2}} \biggl(\pi(\psi_{1})\otimes \pi(\psi_{2})\biggr)
 $$
 $$
 = \biggl(Jac^\theta_{\rho\vert\,\vert^{(a-1)/2}, \cdots
 \rho\vert\,\vert^{(a_{-}+1)/2}}\pi(\psi_{1})\biggr) \otimes \pi(\psi_{2})
 $$
 $$
= \biggl(St(\rho,a_{-}) \times St(\rho,a_{-})\times \pi(\psi'_{1}) \biggr)\otimes \pi(\psi_{2}).
 $$Il y a en plus une action de $\theta$ que l'on n'\'ecrit pas.
Ici on quitte les espaces $I_{cusp}$ \`a cause du fait que $St(\rho,a)$ intervient 2 fois. Mais on sait calculer le transfert vers une distribution stable de $U(n_{1}-(a-a_{-})d_{\rho},E/F)\times U(n_{2},E/F)$, cela vaut exactement
$$
St(\rho,a_{-})\times Transfert^{st} (\pi(\psi'_{1})\otimes \pi(\psi_{2})$$
o\`u le $Transfert^{st}$ est le transfert stable r\'eciproque entre $\tilde{G}_{n'_{1}}\times \tilde{G}_{n_{2}}$ et le produit de groupe unitaire correspondant ($n'_{1}$ est la dimension de la repr\'esentation $\psi'_{1}$ et $n_{2}$ celle de $\psi_{2}$). Il faut maintenant faire le transfert endoscopique vers $U(n,E/F)$; ce transfert est compatible \`a l'induction par $St(\rho,a_-)$ et son image est donc $St(\rho,a_-)\times \pi^H(\psi'_{1}\times \psi_{2})$ o\`u $\pi^H(\psi'_{1}\times \psi_{2}$ est le transfert endoscopique de $Transfert^{st} (\pi(\psi'_{1})\otimes \pi(\psi_{2})$; pour calculer cette repr\'esentation virtuelle, on peut de nouveau se placer dans le $I_{cusp}$ convenable et on obtient un \'el\'ement de la forme $f_{\psi'_{1}\times \psi_{2},s'}$, 
o\`u $s'$ est la restriction de $s$ sur l'espace de $\psi'_{1}\times \psi_{2}$.

On fait le m\^eme calcul en partant de $sz_{\rho,a}$; la d\'ecomposition de $\psi_{s}$ est $\psi'_{1}\times \psi'_{2}$, o\`u $\psi'_{1}$ est comme ci-dessus et $\psi'_{2}$ est la somme de $\psi_{2}$ avec les 2 
repr\'esentations $\rho\otimes \sigma_{[a]}$ et $\rho\otimes \sigma_{[a_{-}]}$. Les calculs se font ensuite exactement comme ci-dessus, c'est $\psi'_{1}$ qui joue un r\^ole muet et $\pi(\psi_{2})$ dont on prend le module de Jacquet. On obtient encore
$$
 Jac^\theta_{\rho\vert\,\vert^{(a-1)/2}, \cdots
 \rho\vert\,\vert^{(a_{-}+1)/2}} \biggl(\pi(\psi'_{1})\otimes \pi(\psi'_{2})\biggr)
 $$
 $$
 = \pi(\psi'_{1})\otimes\biggl(Jac^\theta_{\rho\vert\,\vert^{(a-1)/2}, \cdots
 \rho\vert\,\vert^{(a_{-}+1)/2}} \pi(\psi'_{2})\biggr) =
 $$
 $$
 \pi(\psi'_{1}) \otimes  \biggl(St(\rho,a_-) \times St(\rho,a_-)\times \pi(\psi_{2})\biggr).
 $$
La fin est exactement comme ci-dessus, on commence par sortir l'induction puis on peut retravailler dans $I_{cusp}$. On obtient encore comme image $St(\rho,a_-)\times f_{\psi'_{1}\times \psi_{2},s'}$.

Cela prouve le lemme.
\subsection{Traduction des propri\'et\'es des modules de Jacquet\label{traduction}}
Fixons $\psi$ un morphisme de $W_{E}\times SL(2,{\mathbb C})$ dans $GL(n,{\mathbb C})$, 
$\theta$-discret.

On a d\'ej\`a d\'efini les fonctions $d(\tau,s^{.})$ et les signes $\zeta_{s,\rho,a}$. On a vu que pour $s=1$, $d(\tau,1^{.})\neq 0$ pour tout $\tau \in \Pi(\psi)$; c'est le fait que la projection de $tr\, \tau$ sur 
$I_{cusp}^{st}(U(n,E/F)$ est non nulle. La premi\`ere des conjectures est que pour un bon choix de l'extension de $\pi(\psi)$ \`a $\tilde{G}_{n}$, pour tout $\tau$, $d(\tau,1^{.})$ vaut 
$\vert \Pi(\psi)\vert^{-1}$, c'est-\`a-dire $2^{-\ell(\psi)+1}$. En choisissant la m\^eme normalisation pour tous les groupes $\tilde{G}_{m}$ pour $m\leq n$, on normalise ainsi tous les \'el\'ements qui ont servi \`a d\'efinir la matrice $d(\tau,s^{.})$. Les conjectures d'Arthur  sont exprim\'ees en termes de la matrice inverse de la matrice $\bigl(d(\tau,s^{.})\bigr)$. En traduisant cela donne:

$d(\tau,s^{.})2^{\ell(\psi)-1}$ est un signe et l'application $s\mapsto d(\tau,s^{.})2^{\ell(\psi)-1}$
 est un 
caract\`ere de $(Cent_{GL(n,{\mathbb C})}\psi)^\theta$ trivial sur le sous-groupe du centre de $GL(n,
 {\mathbb C})$ form\'e des \'el\'ements $\theta$-invariants.

 \

 Une variante plus faible serait de montrer que pour tout $s,s'$ dans le centralisateur, $d(s,\tau)d(s',\tau)\neq 0$ et le quotient $d(s,\tau)/d(s',\tau)$ ne d\'epend que de $s s'$ ( tout \'el\'ement du centralisateur est de carr\'e 1). Nous nous allons montrer une forme plus faible encore mais qui suffit pour caract\'eriser les repr\'esen\-ta\-tions cuspidales. Soit $(\rho,a)$ comme dans \ref{definitions} avec $a_{-}$ d\'efini. On prend la notation:
 $$
 Jac_{{\rho\vert\,\vert^{(a-1)/2}}, \cdots, \rho\vert\,\vert^{(a_{-}+1)/2}}\tau:=Jac_{\rho\vert\,\vert^{(a_{-}+1)/2}_{E}}\circ \cdots \circ Jac_{\rho\vert\,\vert^{(a-1)/2}_{E}}\tau.
 $$

 \bf Proposition. \sl Soit $\tau\in \Pi(\psi)$ alors $Jac_{{\rho\vert\,\vert^{(a-1)/2}}, \cdots, \rho\vert\,\vert^{(a_{-}+1)/2}}\tau =0$ si et seulement si pour tout $s$ comme ci-dessus $$
 d(s^{.},\tau)=-\zeta_{s,\rho,a}d(s^{.}z_{\rho,a},\tau).
 $$
 \rm

\

On simplifie la notation en rempla\c{c}ant $Jac_{\rho\vert\,\vert^{(a-1)/2}, \cdots, \rho\vert\,\vert^{(a_{-}+1)/2}}$ par $Jac_{(a-1)/2, \cdots, (a_{-}+1)/2}$.
On \'ecrit avec les d\'efinitions, pour tout $\tau\in \Pi(\psi)$ et pour $\rho,a$ comme ci-dessus:
$$
Jac_{{(a-1)/2},\cdots, (a_{-}+1)/2}\tau=\sum_{s^.}d(s^{.},\tau) Jac_{{(a-1)/2, \cdots, (a_{-}+1)/2}} f_{\psi,s}
$$
$$
=\sum_{(s^.,s^.z_{\rho,a})}(d(s^{.},\tau)+\zeta_{s,\rho,a}d(s^{.}z_{\rho,a},\tau)) 
Jac_{{(a-1)/2}, \cdots, (a_{-}+1)/2} f_{\psi,s}.
$$
Il est clair que si la condition de l'\'enonc\'e est v\'erifi\'ee, cette somme est nulle terme \`a terme. C'est la 
r\'eciproque qu'il faut prouver \`a savoir partir de $\tau$ avec 
$Jac_{{(a-1)/2},\cdots, (a_{-}+1)/2}\tau=0$
et en d\'eduire que chaque terme de la somme est nulle. On reprend les calculs plus pr\'ecis des termes $$Jac_{{(a-1)/2}, \cdots, (a_{-}+1)/2} f_{\psi,s}$$ faits dans la preuve du lemme de \ref{definitions}.
On r\'eutilise la filtration qu'Arthur a introduite dans \cite{arthurselecta};  le terme de plus ''haut'' degr\'e est ce qui provient de $I_{cusp}(U(n,E/F))$ et les autres termes de la filtration ont un gradu\'e qui voit les $I_{cusp}(M)$ pour $M$ les Levi de $U(n,E/F)$. On projette (1) sur le plus haut terme de cette filtration; toutes les contributions des $s$ tels que $\rho\otimes \sigma_{[a]}$ et $\rho\otimes \sigma_{[a_{-}]}$ sont dans le m\^eme espace propre pour $s$ disparaissent. Il ne reste que les contributions des $s$ qui 
s\'eparent ces 2 repr\'esentations. On a alors vu que $Jac_{{(a-1)/2}, \cdots, (a_{-}+1)/2} f_{\psi,s}$ sont des 
\'el\'ements de $I_{cusp}(U(n-2d_{\rho},E/F)$ qui proviennent de la distribution stable associ\'ee \`a $\psi'_{1}\times \psi_{2}$ dans le groupe endoscopique \'evident avec l'\'ecriture. Ces \'el\'ements sont donc lin\'eairement 
ind\'ependants et la nullit\'e de (1) assure que pour ces $s$, on a la relation de l'\'enonc\'e. Il ne reste donc plus que les termes faisant intervenir une induction par $St(\rho,a_{-})$ avec un \'el\'ement $f_{\psi'_{1}\times \psi_{2},s'}$. Ces termes sont aussi lin\'eairement ind\'ependants et on obtient la relation
annonc\'ee dans l'\'enonc\'e.
\subsection{Remarques}
On reprend les notations $\psi$, $(\rho,a,a_{-})$ du paragraphe pr\'ec\'edent.

\

\bf Remarque.\rm $\{\tau \in \Pi(\psi); Jac_{\rho\vert\,\vert^{(a-1)/2}, \cdots, \rho\vert\,\vert^{(a_{-}+1)/2}}\tau=0\}$ est de cardinal $\vert \Pi(\psi)\vert/2$ c'est-\`a-dire $2^{\ell(\psi)-2}$.

\

Le cardinal de l'ensemble de l'\'enonc\'e de la remarque est \'egal au rang de la matrice rectangulaire 
$d(\tau,s^.)$, o\`u $\tau$ parcourt cet ensemble et $s^.$ parcourt $(Cent_{GL(n,{\mathbb C})}(\psi)/{\mathbb C}^*)^\theta$ o\`u ${\mathbb C}^*$ est le centre de $GL(n,{\mathbb C})$. Et le lemme 
pr\'ec\'edent dit que cette matrice est de rang 1/2 la taille de la matrice carr\'e $d(\tau,s^.)$ o\`u ici $\tau$ parcourt $\Pi(\psi)$. C'est le r\'esultat annonc\'e.

\subsection{Liens avec les conjectures d'Arthur}
\subsubsection{Pr\'eliminaires \label{preliminaires}}
C'est la matrice inverse de la matrice $d(\tau,s^{.})$ de \ref{definitions} qui intervient plus 
spontan\'ement. On note $c(s^{.},\tau)$ cette matrice inverse; elle est d\'efinie par le fait que pour $s$ un \'el\'ement $\theta$-invariant du centralisateur de $\psi$, on note $\psi_{s}:=\psi_{1}\times \psi_{2}$ la 
d\'ecomposition de 
$\psi$ telle que $\psi_{1}$ soit \`a valeurs dans l'espace propre pour la valeur propre $-1$ de $s$ et $\psi_{2}$ dans l'espace propre pour la valeur propre $+1$. On note $H_{s}$ la donn\'ee endoscopique 
d\'etermin\'ee par $s$. On consid\`ere la repr\'esentation $\pi(\psi_{s})$ et un prolongement de cette 
repr\'esentation au groupe tordu par $\theta$; cela d\'etermine une distribution stable sur $H_{s}$, 
not\'ee $\Pi^{st}_{H_{s}}(\psi_{s})$, combinaison lin\'eaire de caract\`eres de s\'eries discr\`etes de 
$H_{s}$. Il faut regarder simultan\'ement $\psi_{s}$ et $\psi_{zs}$ o\`u $z$ est la matrice diagonale n'ayant que la valeur propre $-1$. Le choix de $\tilde{\pi}(\psi_{s})$ induit un choix de $\tilde{\pi}(\psi_{zs})$ de telle sorte que $\Pi^{st}_{H_{s}}(\psi_{s})+\Pi^{st}_{H_{zs}}(\psi_{zs})$ soit invariant par le groupe des automorphismes venant de $U(n,E/F)$. On note $\Pi^{s^.-st}(\psi)$ l'image de cet \'el\'ement dans $I_{cusp}(U(n,E/F))$ et on voit cet \'el\'ement comme une combinaison lin\'eaire de caract\`etres de repr\'esentations dans $\Pi(\psi)$. Alors, on a par d\'efinition:
$$
\Pi^{s^.-st}(\psi)=\sum_{\tau\in \Pi(\psi)}c(s^.,\tau) \tau.
$$On reprend la notation $\zeta_{s,\rho,a}$ de \ref{traduction} et on a:

\

\bf Corollaire. \sl Soit $\tau\in \Pi(\psi)$ et $(\rho,a,a_{-})$ comme en \ref{definitions}. Alors, 
$$
Jac_{\rho\vert\,\vert^{(a-1)/2}, \cdots, \rho\vert\,\vert^{(a_{-}+1)/2}}\tau \neq 0\eqno(1)
$$
si et seulement si pour tout $s^.$, $c(s^.,\tau)=\zeta_{s,\rho,a}c(s^.z_{\rho,a},\tau).
$
\rm

\

Par d\'efinition pour tout $\tau'\in \Pi(\psi)$ et tout $\tau\in \Pi(\psi)$, on a:
$$
\sum_{s^.}d(\tau',s^.)c(s^.,\tau)=\delta_{\tau',\tau}\eqno(2)
$$
o\`u $\delta_{\tau',\tau}$ vaut $0$ si $\tau'\neq \tau$ et 1 si $\tau'=\tau$. Ainsi la condition de l'\'enonc\'e est \'equivalente \`a ce que (2) soit $0$ pour tout $\tau'$ tel que $Jac_{\rho\vert\,\vert^{(a-1)/2}, \cdots, \rho\vert\,\vert^{(a_{-}+1)/2}}\tau'=0$. Fixons un tel $\tau'$.
On regroupe $s^.$ et $s^.z_{\rho,a}$ et en appliquant \ref{traduction}, la condition (1) est \'equivalente \`a, pour tout $\tau'$ comme ci-dessus:
$$
\sum_{s^.,s^.z_{\rho,a}}d(\tau',s^.)(c(s^.,\tau)+\zeta_{\rho,a}c(s^.z_{\rho,a},\tau))=0
$$
Comme la matrice $d(\tau',s^.)$ o\`u $\tau'$ est comme ci-dessus et $s^.$ parcourt un ensemble de 
repr\'esentants modulo la multiplication par $z_{\rho,a}$ est de rang maximum, on obtient le corollaire.

\

Pour ne pas compliquer inutilement les notations, on note $<\psi_{s},\psi_{s}>$ la norme elliptique de la repr\'esentation virtuelle $\Pi^{s^.-st}(\psi)$, c'est \`a dire la norme dans $I_{cusp}(U(n,E/F)$. On a:

\

\bf Lemme. \sl Pour tout $s,\tau$ comme ci-dessus, $c(s^.,\tau)=d(\tau,s^.)<\psi_{s},\psi_{s}>$.\rm

\

\noindent
Soient $\tau,\tau'\in \Pi(\psi)$; on \'ecrit $<tr\, \tau,tr\,\tau'>=\delta_{\tau,\tau'}$ o\`u $<\, , \, >$ est le produit scalaire dans $I_{cusp}(U(n,E/F)$ et cela donne:
$$
\sum_{s^.}d(\tau,s^.)d(\tau',s^.)<\psi_{s},\psi_{s}>=\delta_{\tau,\tau'}.
$$
Le lemme r\'esulte donc de la d\'efinition de $c(s^.,\tau)$ comme inverse de la matrice $d(\tau,s^.)$.

\

\bf Corollaire. \sl
Soit $\tau\in \Pi(\psi)$ et $(\rho,a,a_{-})$ comme en \ref{definitions}. Alors, 
$$
Jac_{\rho\vert\,\vert^{(a-1)/2}, \cdots, \rho\vert\,\vert^{(a_{-}+1)/2}}\tau= 0\eqno(2)
$$
si et seulement si pour tout $s^.$, $c(s^.,\tau)=-\zeta_{s,\rho,a}c(s^.z_{\rho,a},\tau)<\psi_{sz_{a}},\psi_{sz_{a}}><\psi_{s},\psi_{s}>^{-1}.
$
\rm

\

C'est un corollaire imm\'ediat de \ref{traduction}

\subsubsection{Normalisation}

Jusqu'\`a pr\'esent, nous avons \'evit\'e de normaliser l'action de $\theta$ sur $\pi(\psi)$, puisque cela n'\'etait pas utile. On pourrait continuer ainsi mais cela complique singuli\`erement les d\'efinitions; on va donc montrer que la normalisation la plus populaire (\`a l'aide de mod\`ele de Whittaker) permet de faire dispara\^{\i}tre le signe $\zeta_{s,\rho,a}$. On fixe donc un caract\`ere additif de $E^*$ et un caract\`ere $\theta$ invariant non d\'eg\'en\'er\'e du groupe unipotent sup\'erieur de $GL(n,E)$; c'est le caract\`ere usuel puisque $\theta$ respecte un \'epinglage. Pour faire ceci, il suffit de supposer que la classe de conjugaison de $\psi$ est invariante sous l'action de $\theta$ (vu dualement) et il n'est pas utile de supposer que $\psi$ est $\theta$ discret.

On a \'etudi\'e les diff\'erentes normalisations en \cite{transfert} et il r\'esulte de ces constructions que dans le cas des repr\'esentations temp\'er\'ees tous les choix raisonables donnent le m\^eme r\'esultat. C'est donc ce qui se produit ici aussi et va permettre de prouver que pour le choix ci-dessus, les signes $\zeta_{s,\rho,a}$ sont tous $+1$.

Rappelons une construction simple: soit $\sigma$ une repr\'esentation d'un groupe lin\'eaire $GL(n',E)$ avec une action de $\theta$, not\'ee $\theta_{\sigma}$. Et soit $\lambda$ une repr\'esentation d'un groupe $GL(\ell,E)$; on note $\theta(\lambda)$ l'image de $\lambda$ par $\theta$. La repr\'esentation induite $\lambda\times \sigma\times \theta(\lambda)$ a une action naturelle de $\theta$ provenant de $\theta_{\sigma}$. En effet, on fixe $A_{\lambda}$ un homomorphisme de l'espace o\`u $\lambda$ se 
r\'ealise dans l'espace o\`u $\theta(\lambda)$ se r\'ealise v\'erifiant, pour tout $v$ dans l'espace de $\lambda$ et tout $g\in GL(\ell,F)$
$$
A_{\lambda}(\lambda(\theta(g)).v)=\theta(\lambda)(g)A_{\lambda}(v).
$$
Soit $f\in \lambda\times \sigma\times \theta(\lambda)$, c'est-\`a-dire une fonction sur $GL(n'+2\ell,E)$ \`a valeurs dans l'espace de la repr\'esentation $\lambda\otimes \sigma\otimes \theta(\lambda)$. On pose:
$$
\theta(f):= (g\in GL(n'+2\ell,E) \mapsto  (A_{\lambda}^{-1}\otimes \theta_{\sigma} \otimes A_{\lambda})\circ inv \, f(\theta(g))
$$
o\`u $inv$ est l'application qui envoie $\lambda\otimes \sigma\otimes \theta(\lambda)$ dans $\theta(\lambda)\otimes \sigma \otimes \lambda$ en \'echangeant les facteurs extr\^emes. Cette action de $\theta$ ne d\'epend que de $\theta_{\sigma}$ et non pas du choix de $A_{\lambda}$.

Soit ${\cal E}$ une collection de couples $(\rho,a)$ avec $\rho$ une repr\'esentation cuspidale 
irr\'eductible $\theta$ invariante et $a$ un entier. On fixe $(\rho',a')\in {\cal E}$ avec $a'\geq 2$ et on pose ${\cal E}'$ l'ensemble qui se d\'eduit de ${\cal E}$ en rempla\c{c}ant $(\rho',a')$ par $(\rho',a'-2)$. On pose:
$$
\pi:=\times_{(\rho,a)\in {\cal E}}St(\rho,a); \pi':=\times_{(\rho,a)\in {\cal E}'}St(\rho,a).
$$
Ces 2 repr\'esentations sont munies d'une action de $\theta$ fixant le mod\`ele de Whittaker; on note $\theta^W_{\pi}$ et $\theta^W_{\pi'}$ ces actions. On sait aussi que $\pi$ est un sous-module irr\'eductible de l'induite
$$
\sigma:=\rho'\vert\,\vert^{(a'-1)/2}_{E}\times \pi' \times \rho'\vert\,\vert^{-(a'-1)/2}. \eqno(1)
$$
Et $\pi$ est l'unique sous-quotient de cette induite ayant un mod\`ele de Whittaker, ce qui prouve 
imm\'e\-dia\-tement le fait que $\pi$ intervient dans cette induite avec multiplicit\'e 1.
Comme on a suppos\'e que $\rho'$ est $\theta$ invariante, l'induite (1) a une action de $\theta$ qui provient de $\theta^W_{\pi'}$. On note $\theta$ cette action.

\

\bf Lemme. \sl l'action $\theta^W_{\pi}$ est la restriction de $\theta$ \`a $\pi$ pour l'inclusion de $\pi$ dans (1).\rm

\

On note $V'$ l'espace de $\rho'\vert\,\vert^{(a'-1)/2}_{E}\otimes \pi' \otimes \rho'\vert\,\vert^{-(a'-1)/2}$ et $\ell^W_{V'}$ une fonctionnelle de Whittaker sur $V'$; cela \'etend la fonctionnelle de Whittaker sur $\pi'$ de fa\c{c}on invariante sous $inv$ utilis\'e dans les constructions ci-dessus. On construit une fonctionnelle de Whittaker sur $V$, l'espace de $\sigma$, de fa\c{c}on standard (on 
 rappelera sa construction) $\ell^W_{V}$. On note $n$ l'entier tel que $\pi$ soit une repr\'esentation de $GL(n,E)$ et $N$ le groupe des matrices unipotentes sup\'erieures de $GL(n,E)$; on note $\chi$ le caract\`ere de $N$, $\theta$-invariant qui sert \`a d\'efinir le mod\`ele de Whittaker. D'abord $\theta^W_{\pi}$ est l'action de $\theta$ qui induit l'identit\'e sur l'espace vectoriel de dimension 1, $V/N_{\chi}V$ (module de Jacquet tordu). Pour construire $\ell^W_{V}$ on identifie $V$ \`a un espace de fonction de $GL(n,E)$ \`a valeurs dans $V'$, covariante pour l'action du sous-groupe parabolique sup\'erieur de Levi isomophe \`a $M:=GL(d_{\lambda},E) \times GL(n',E)\times GL(d_{\lambda},E)$ (avec des notations \'evidentes). On note $w$ l'\'el\'ement du groupe de Weyl de $GL(n,E)$ de longueur minimale dans sa double classe modulo le groupe de Weyl de $M$ qui rend n\'egatif toutes les racines simples hors de $M$.

On note $N_{P}$ le radical unipotent du parabolique qui sert \`a l'induction. Et on pose pour tout $f\in V$:
$$
\ell^W_{V}(f):=\int_{N_{P}}\ell^W_{V'}f(wn)\chi(n) \, dn.
$$
L'int\'egrale se fait sur un compact et est donc une somme finie qui se factorise par $V/N_{\chi}V$. Pour 
v\'erifier que l'action de $\theta$ sur (1) prolongeant canoniquement $\theta^W_{\pi'}$ induit aussi l'action triviale sur $V/N_{\chi}V$, il suffit de calculer explicitement $\ell^W_{V}(\theta.f)$; le r\'esultat est la m\^eme int\'egrale mais avec $w$ remplac\'e par $\theta(w)$. On v\'erifie que $\theta(w)=zw$ avec $z$ un \'el\'ement du centre de $M$ de la forme $z_{1}\times 1 \times z_{1}$ avec $z_{1}$ valant 
$\pm Id_{d_{\lambda}}$. La repr\'esentation induisante a un caract\`ere central trivial sur $z$ d'o\`u l'invariance
$$
\ell^W_{V}(f)=\ell^W_{V}(\theta(f)).
$$
Cela assure que $\theta(f)$ et $f$ ont m\^eme image dans $V/N_{\chi}V$ comme cherch\'e.

\

\bf Corollaire. \sl Toutes les actions de $\theta$ \'etant fix\'ees par la normalisation de Whittaker, les signes $\zeta_{s,\rho,a}$ de \ref{traduction} sont tous \'egaux \`a $+1$.\rm

\

Il faut d'abord se ramener au cas o\`u $s=1$; cela r\'esulte du fait que le transfert commute \`a la prise du module de Jacquet (cf. \ref{moduledejacquet}. Pour $s=1$, avec $(\rho,a)$ comme dans l'\'enonc\'e on vient de montrer que l'inclusion:
$$
\times_{(\rho',a')\in {\cal E}}St(\rho,a) \hookrightarrow \rho\vert\,\vert^{(a-1)/2}_{E} \times \biggl(\times_{(\rho',a')\in {\cal E}-\{(\rho,a)\}}St(\rho',a') \times St(\rho,a-2)\biggl)\times \rho\vert\,\vert^{-(a-1)/2}_{E}\eqno(2)
$$
est compatible aux actions de $\theta$ quand sur le membre de droite on met l'action prolongeant canoniquement l'action de $\theta$ sur la repr\'esentation entre parenth\`ese. Ce prolongement canonique est compatible au module de Jacquet, par d\'efinition puisque
$$Jac^\theta_{\rho\vert\,\vert^{(a-1)/2}_{E}}(2)=(\times_{(\rho',a')\in {\cal E}-\{(\rho,a)\}}St(\rho',a') \times St(\rho,a-2)$$ 
s'obtient en \'evaluant \`a l'origine. Si $a_{-}$ (notation de \ref{traduction}) vaut $a-2$, on a 
imm\'ediatement $\zeta_{1,\rho,a}=+1$ et sinon on obtient le corollaire en r\'eit\'erant cette construction.

\subsubsection{Caract\`ere de $A(\psi)$\label{caracterebis}}
On fixe $\psi$ un morphisme $\theta$-discret de $W_{E}\times SL(2,{\mathbb C})$ dans 
$GL(n,{\mathbb C})$. D'o\`u un paquet de repr\'esenta\-tions $\Pi(\psi)$ de $U(n,E/F)$. On a introduit 
$A(\psi)$ un sous-groupe du centralisateur de $\psi$ en \ref{preliminaires}. Pour comprendre ce qui se passe, on va tout de suite construire un sous-groupe $S_{0}$ du centralisateur de $\psi$ de telle sorte que:
$$S_{0}\cap A(\psi)=\{1\}, \qquad
(Cent_{GL(n,{\mathbb C})}\psi)^\theta=S_{0}A(\psi). \eqno(1)
$$
On \'ecrit $\psi$ comme une somme de repr\'esentations irr\'eductibles $\rho\otimes \sigma_{[a]}$ de $W_{E}\times SL(2,{\mathbb C})$ et on note ${\cal E}$ l'ensemble des couples $(\rho,a)$ intervenant. On sait que le centralisateur de $\psi$ dans $GL(n,{\mathbb C})$ est le produit des groupes 
${\mathbb C}^*$ agissant scalairement sur chaque espace des sous-repr\'esentations $\rho\otimes \sigma_{[a]}$; quand on consid\`ere les \'el\'ements $\theta$-invariants de ce groupe, il reste ${\pm 1}$. On consid\`ere ${\cal E}_{imp}$ l'ensemble des couples $(\rho,a)\in {\cal E}$ tel que $a$ soit impair. On rappelle que cette condition est \'equivalente \`a ce que $L(\rho,r_{A}',s)$ n'ait pas de p\^ole en $s=0$. On remarque que $A(\psi)$ contient tous les facteurs ${\pm 1}$ relatif aux \'el\'ements de ${\cal E}-{\cal E}_{imp}$. Donc en particulier si ${\cal E}_{imp}=\emptyset $, (1) est r\'ealis\'e avec $S_{0}=\emptyset$ et $Z$ est alors inclus dans $A(\psi)$. Supposons donc que ${\cal E}_{imp}\neq \emptyset$. On fixe $\rho$ et un nombre impair d'entier ${A}_{\rho}$ tel que $(\rho,a)\in {\cal E}_{imp}$ pour tout $a\in A_{\rho}$. On note $z_{\rho}$ l'\'el\'ement du centralisateur de $\psi$ produit des \'el\'ements $-1$ relatifs \`a ces repr\'esentations $\rho\otimes\sigma_{[a]}$ pour $a$ parcourant $A_{\rho}$. On note alors $S_{0}$ le groupe engendr\'e par ces \'el\'ements $z_{\rho}$. Evidemment la d\'efinition n'est pas du tout canonique puisqu'elle d\'epend du choix de $A_{\rho}$ mais (1) est clair. Un choix qui en vaut bien un autre est de prendre $A_{\rho}$ r\'eduit \`a un \'el\'ement et pr\'ecis\'ement \`a l'\'el\'ement le plus petit possible.

\

\bf Th\'eor\`eme. \sl Pour tout $\tau \in \Pi(\psi)$, il existe un unique caract\`ere, $\epsilon_{\tau}$ de $A(\psi)$ tel que pour tout $s\in (Cent_{GL(n,{\mathbb C})}\psi)^\theta$ et tout $a\in A(\psi)$, il existe un 
r\'eel strictement positif $\alpha$ avec l'\'egalit\'e:
$$
c(s^.a,\tau)=\epsilon_{\tau}(a) c(s^.,\tau) \alpha.
$$
\rm

\

Montrons d'abord l'unicit\'e de $\tau$: on fixe $s=1$, on sait alors que $c(1^.,\tau)\neq 0$ et $\epsilon_{\tau}(a)$ est alors le signe de $c(a^.,\tau)/c(1^.,\tau)$ qui est calcul\'e en fonction des modules de Jacquet de $\tau$ gr\^ace aux 2 propri\'et\'es d\'emontr\'ee en \ref{preliminaires}. On peut donc remplacer $1^.$ dans cette d\'efinition par n'importe quel \'el\'ement de $A(\psi)$. Remarquons maintenant que cette application signe d\'efinit donc bien un caract\`ere. En appliquant ces m\^emes 
r\'ef\'erences,  o\`u on peut prendre $s$ quelconque, on obtient le th\'eor\`eme. On conjecture que $c(s^.,\tau)=\pm 1$ et que $\alpha=1$ mais ce n'est pas prouv\'e ici. Toutefois, on peut calculer $\alpha$ comme quotient de normes (que l'on ne sait pas calculer).

\subsubsection{Classification des repr\'esentations cuspidales\label{classificationcuspidales}}
On reprend les notations $\psi$, $A(\psi)$ et $\epsilon_{\tau}$ de \ref{caracterebis}. On note $\epsilon_{alt}$ l'unique caract\`ere de $A(\psi)$ qui vaut $-1$ sur  tout g\'en\'erateur $z_{\rho,a}$ d\'efini  en 
\ref{traduction}. On a d\'efini le fait que $\psi$ soit sans trou dans l'introduction; on a vu que si $\tau\in \Pi(\psi)$ le support cuspidal \'etendu de $\tau$ se calcule avec la d\'ecomposition en repr\'esentations irr\'eductibles de $\psi$. On a donn\'e la forme du support cuspidal \'etendu d'une repr\'esentation cuspidal en \ref{supportcuspidalcuspidal} et il en r\'esulte que les seuls morphismes $\psi$ tels que $\Pi(\psi)$ contienne une repr\'esentation cuspidal sont les $\psi$ qui sont sans trou.

\

\bf Th\'eor\`eme. \sl On suppose que $\psi$ est sans trou. La repr\'esentation $\tau$ est cuspidale si et seulement si 
$\epsilon_{\tau}=\epsilon_{alt}$. Le nombre de repr\'esentations cuspidales dans $\Pi(\psi)$ est 
le cardinal de l'ensemble des caract\`eres de $(Cent_{GL(n,{\mathbb C})}\psi)^\theta/\{Id,-Id\}$ 
dont la restriction \`a $A(\psi)$ vaut $\epsilon_{alt}$.\rm

\

La premi\`ere partie du th\'eor\`eme r\'esulte de \ref{preliminaire} montrant que $\epsilon_{\tau}=\epsilon_{alt}$ si et seulement si tous les modules de Jacquet de $\tau$ sont nuls. Montrons la 2e partie: le premier point \`a remarquer est que si $Z=\{Id,-Id\}$ est inclus dans $A(\psi)$, la restriction de $\epsilon_{\tau}$ \`a $Z$ est n\'ecessairement l'identit\'e. Il y a 2 cas \`a distinguer, le premier cas est celui o\`u $Z$ est inclus dans $A(\psi)$ et o\`u $\epsilon_{alt}$ n'est pas de restriction triviale  \`a $A(\psi)$. Dans ce cas, $\Pi(\psi)$ ne contient pas de repr\'esentations cuspidales et la fin du th\'eor\`eme est claire.
Dans le cas oppos\'e, on note encore $\epsilon_{alt}$ le caract\`ere du groupe $ZA(\psi)$ qui vaut $\epsilon_{alt}$ sur $A(\psi)$ et est trivial sur $Z$. On fixe $S'_{0}$ un sous-ensemble de $S_{0}$ tel que $S_{0}A(\psi)=S'_{0}A(\psi)Z$. Pour tout $\tau\in \Pi(\psi)$ on note encore $\epsilon_{\tau}$ le caract\`ere de $ZA(\psi)$ dont la restriction \`a $A(\psi)$ est $\epsilon_{\tau}$ et qui est trivial sur $Z$.
On consid\`ere l'application qui \`a $\tau\in \Pi(\psi)$ associe l'\'el\'ement $(c(s'_{0},\tau)\in {\mathbb C}^{\vert S'_{0}\vert}),\epsilon_{\tau}$ c'est \`a dire un \'el\'ement de 
${\mathbb C}^{\vert S'_{0}\vert}$ et un caract\`ere de $ZA(\psi)$ trivial sur $Z$. Le fait que la matrice est de rang $c(s^.,\tau)$ est inversible se traduit par le fait que pour tout caract\`ere $\epsilon_{0}$ de 
$ZA(\psi)$ trivial sur $Z$, la matrice
$\{c(s^._{0},\tau); s_{0}\in S'_{0}, \epsilon_{\tau}=\epsilon_{0}\}$ est de rang $\vert S'_{0}\vert$. Ceci est a fortiori vrai pour $\epsilon_{0}=\epsilon_{alt}$ et donne le r\'esultat cherch\'e.

\end{document}